\documentclass{article}

\usepackage[T1]{fontenc}

\usepackage{amsfonts}
\usepackage{amsmath}
\usepackage{cite}
\usepackage{graphicx}
\usepackage{lipsum}
\usepackage{caption}
\usepackage{subcaption}
\usepackage{algorithm}
\usepackage{algpseudocode}
\usepackage{multirow}
\usepackage[a4paper, total={6in, 8in}]{geometry}
\usepackage[table,xcdraw]{xcolor}
\usepackage{enumitem}
\usepackage{tabularx}
\usepackage{pifont}
\newcommand{\cmark}{\ding{51}}%
\usepackage{tablefootnote}
\usepackage{hyperref}
\usepackage{longtable}
\usepackage{pgfplots} 
\usepackage{svg}
\usepackage{authblk} 
\usepackage{tikz}
\usetikzlibrary{arrows.meta, positioning}
\usepackage{booktabs}
\usepackage{multirow}
\usepackage{array}
\usepackage{siunitx}
\usepackage{multicol} 

\usepackage{booktabs}
\usepackage{multirow}
\usepackage{threeparttable}
\usepackage{geometry}

\usepackage{natbib} 
\usepackage{hyperref} 

\usepackage{tabularx}
\usepackage{ragged2e} 

\allowdisplaybreaks

\title{RoRo stowage planning: a systematic study of mathematical formulations.}
\author[1,3]{Oliver Rise Thomsen}
\author[1,2,3]{Magnus Vilhelm Barrett Levinsen}
\author[1]{Dario Pacino}
\author[2]{Line Blander Reinhardt}

\affil[1]{Department of Technology, Management and Economics, Technical University of Denmark}
\affil[2]{Roskilde University}
\affil[3]{DFDS A/S}

\begin{document}

\maketitle

\begin{abstract}

Stowage planning for Roll‑On/Roll‑Off (RoRo) vessels can be represented with different mathematical formulations for efficient stowage plans: hybrid‑slot, grid‑based, and lane‑based representations. In this work, we implement and compare three formulations equally across a diverse set of instances. We extend these three representations to handle heterogeneous cargo in terms of size and weight, and to address advanced stability requirements, such as centre of gravity limits, stress force limits, and ballast‑water adjustments. We evaluate the models in three configurations: cargo placement, simple stability with centres of gravity limits and ballast adjustments, and full stability, which adds shear‑force and bending‑moment limits. The tests are conducted on three different vessel sizes and solved with a mixed‑integer programming solver. The hybrid‑slot and grid‑based models both perform well finding solutions for 447 out of 486 instances across all configurations; however, the hybrid‑slot model placed 8\% more cargo than the grid model. By contrast, the lane model only finds solutions for 142 instances and placed 15\% less cargo than the grid‑based model. The solutions generated by the hybrid‑slot model are validated using a commercial vessel‑stability software to verify our stability estimations, showing issues in only 24 of 162 instances, which are easily avoided by tightening the stability constraints. Finally, we address the placement of hazardous and refrigerated cargo in an efficient manner for the hybrid‑slot model by employing lazy constraints which is shown to solve realistic sized instances.

\end{abstract}

\newpage

\section{Introduction}\label{sec:introduction}
Transport of wheeled cargo by sea, such as truck trailers, cars, and agricultural equipment, is known as Roll-on/Roll-off (RoRo) shipping. Accounting for 12.7\% of the gross cargo tonnage at major European ports, RoRo shipping plays a vital role in both short-sea and deep-sea transportation. This significance is further affirmed by the fact that in Ireland and Estonia, RoRo represents over 30.0\% of the total transported tonnage \cite{eurostat_maritime_2025}.

Most commonly, RoRo vessels are deployed on short-sea shipping routes and are an integral part of intermodal transport chains, where the main cargo consists of trailers or specialised wheeled platforms of bulk equipment. The deployment of ocean sailing RoRo vessels is often dedicated to transporting cars and heavy machinery from manufacturers to consumers.
Unlike container shipping, RoRo vessels use ramps that allow cargo to be rolled on/off directly, eliminating the need for large quay cranes. The wheeled cargo is driven directly onto the vessel either by itself or pulled using tugs (terminal tractors).

Research shows that higher automation and better planning at RoRo terminals can reduce turnaround times and improve vessel utilisation \cite{ovstebo_routing_2011,hansen_adaptive_2020, puisa_optimal_2021, lalla-ruiz_stowage_2020,jia_optimal_2022, marquardt_continuous_2026}. The focus of this paper is specific to stowage planning.
RoRo stowage planning concerns the identification of where cargo must be positioned on a vessel. Unlike container shipping, where three indices (bay, stack, and tier) identify the position of a 40-foot container \cite{van_twiller_literature_2023}, RoRo vessels have open deck spaces where cargo can be placed virtually anywhere, similar to two-dimensional packing problems \cite{hansen_2d-packing_2016}. The decks can only be accessed by ramps. Lanes, running the entire length of the ship, are painted on the deck floor and are often used to assist in planning \cite{ovstebo_optimization_2011}. The cargo can have different shapes and dimensions, complicating stowage and making it challenging to achieve a space-efficient arrangement. Cargo might also have special requirements that constrain its position within the vessel (e.g., access to electrical outlets) or with respect to other cargo (e.g., dangerous goods). Moreover, the positioning of cargo on the vessel must ensure that it remains seaworthy and satisfies all stability and stress limits. 

The overall aim of RoRo ship stowage planning is the minimisation of costs while ensuring safe operations, typically measured in terms of cargo-handling efficiency \cite{daduna_dynamic_2023} or bunker reduction (e.g. the minimisation of ballast water) \cite{lalla-ruiz_stowage_2020}. 

Compared with the Container Stowage Planning Problem (CSPP) \cite{van_twiller_literature_2023}, the literature on RoRo Stowage Planning remains relatively sparse. This paper aims to analyse and compare the three principal mathematical representations used for RoRo stowage planning, the Lane, Grid, and Hybrid-Slot approaches, with particular emphasis on their ability to allocate individual cargo items and their performance under seaworthiness constraints. The main contributions of this paper are as follows:

\begin{itemize}
    \item The extension of the state-of-the-art mathematical formulations to include a common set of industrially applicable sea-worthiness conditions, including shear forces, bending moments, and hydrostatic data.
    \item Extension and reformulation of the lane model to include general shaped vessel layouts.
    \item A comparative evaluation of the modelling approaches on real-life vessel instances.
    \item Modelling of refrigerated and hazardous cargo in the Hybrid-Slot model.
    \item Quality assessment of the linearization of the stability constraints against the precise calculations of a commercial stability software.
\end{itemize}

The remainder of this paper is structured as follows: \autoref{sec:literature} introduces the literature on the RoRo Ship Stowage Problem (RSSP), providing an overview of the existing research. \autoref{sec:formulation} presents the RSSP, while the various representations, and the new stability and stress forces are modelled in \autoref{sec:math}. \autoref{sec:results} offers a comprehensive analysis of the three modelling approaches. In \autoref{sec:extension}, we demonstrate how the hybrid-slot model can be extended to account for hazardous and refrigerated cargo, and test the model's seaworthiness against a commercial stability software. 
Finally, the conclusions are drawn in \autoref{sec:conclusion}.

\section{Related Work} \label{sec:literature}
The research field of RoRo stowage planning optimisation is relatively new. As mentioned earlier, fewer studies are available compared to its containerised equivalent, which was most recently reviewed in \cite{van_twiller_literature_2023}. The key difference between container and RoRo stowage is that RoRo vessels transport wheeled cargo with varying dimensions, to be stowed in a more or less open floor plan (in contrast to the cellular structure of container vessels). The open floor plan makes the modelling of the cargo position more challenging. The first paper that introduced the RoRo Ship Stowage Problem (RSSP) appeared in 2011 \cite{ovstebo_optimization_2011}. The problem has since received increasing interest. Multiple papers have emerged that have studied different aspects of the problem, introduced various model formulations, and implemented new solution methods. \autoref{tab:literature} shows these publications and their main differences.

\newcolumntype{Y}{>{\centering\arraybackslash}X}

\begin{table}[!thb]
\centering
\scriptsize
    \begin{tabularx}{\textwidth}{|l|Y|Y|Y|Y|Y|Y|}
    \hline
     & \textbf{2011} \par \cite{ovstebo_routing_2011}, \cite{ovstebo_optimization_2011} 
     & \textbf{2016-2020} \par \cite{hansen_2d-packing_2016}, \cite{hansen_adaptive_2020} 
     & \textbf{2018} \par \cite{puisa_optimal_2021} 
     & \textbf{2020-2021} \par \cite{lalla-ruiz_stowage_2020}, \cite{jia_step-wise_2021} 
     & \textbf{2024} \par \cite{daduna_dynamic_2023} 
     & \textbf{This paper} \\ \hline
\textbf{Representation}          & Lane                                                                          & Grid                                                                              & Hybrid-Slot                                                           & Slot                                                                            & Grid                                                                       & All                                                                \\ \hline
\textbf{Variable cargo sizes}    & \cmark                                                                        & \cmark                                                                               & \cmark                                                                   & -                                                                              & \cmark                                                                        & \cmark                                                                \\ \hline
\textbf{Free cargo positioning}  & (\cmark$^1$)                                                                          & \cmark                                                                               & -                                                                    & -                                                                              & \cmark                                                                        & (\cmark$^2$)                                                            \\ \hline
\textbf{Multi-port voyages}      & \cmark                                                                           & \cmark                                                                               & \cmark                                                                   & -                                                                              & -                                                                         &  -                                                                 \\ \hline
\textbf{Multiple deck}           & \cmark                                                                           & -                                                                                & \cmark                                                                   & \cmark                                                                             & \cmark                                                                        & \cmark                                                                \\ \hline
\textbf{Uncertain cargo arrival} & -                                                                            & -                                                                                & -                                                                    & -                                                                              & \cmark                                                                        & -                                                                 \\ \hline
\textbf{Stability}               & \cmark                                                                           & -                                                                                & \cmark                                                                   & \cmark                                                                             & -                                                                         & \cmark                                                                \\ \hline
\textbf{Ballast Tanks}           & -                                                                            & -                                                                                & -                                                                    & \cmark                                                                             & -                                                                         & \cmark                                                                \\ \hline
\textbf{Hazardous Cargo}         & -                                                                            & -                                                                                & -                                                                    & \cmark                                                                             & -                                                                         & \cmark                                                                \\ \hline
\end{tabularx}
\caption{The figure illustrates the RSSP literature, highlighting the primary objectives and constraints addressed by each paper. 1. The lane representation only has free cargo positioning the the longitudinal direction. 2. This papers models all three representations, thus it depends on the model.}
\label{tab:literature}
\end{table}

A significant difference in the published research is the representation of the deck area of the vessel. This representation determines how cargo positions are defined. The lane-based approach was introduced in \cite{ovstebo_optimization_2011}, where the layout of a vessel deck is split into a set of lanes in which cargo can be placed. This layout representation is inspired by the industry and the use of painted lanes within the vessel decks. However, the model allowing lanes to have a width dependent on the width of the cargo placed within the lane, by \cite{ovstebo_routing_2011}, limits the possible deck designs to a rectangle. The grid representation was introduced in \cite{hansen_2d-packing_2016}, where the deck area is overlaid by a grid (or a mesh of points). Each node in the grid represents a possible position of a cargo. The grid representation allows for better area utilisation if the grid resolution is granular; however, increased granularity of the grid will also increase the complexity of the model.
The hybrid-slot deck representation divides the deck area into slots (discrete cargo positions) for each cargo type, ensuring that slots in the same area do not overlap. The hybrid-slot approach was presented in \cite{puisa_optimal_2021}, and a simplification of this approach (the slot representation), where only homogeneous cargo is modelled, was presented in \cite{lalla-ruiz_stowage_2020}. State-of-the-art mathematical formulations, based on the described deck layout representation, can only efficiently be used to solve small-scale instances. Only the slot and the hybrid representation have been shown to find optimal solutions to industrial-sized problems utilising an MIP solver. However, since all published studies focus on different variations of the RoRo stowage planning problem, a systematic comparison of these models is not available. 

Most of the published research categorises cargo into mandatory and optional, aiming to maximise the revenue from optional cargo \cite{ovstebo_optimization_2011, hansen_2d-packing_2016, puisa_optimal_2021, daduna_dynamic_2023}. Moreover, the literature agrees that stowage plans should be able to handle heterogeneous cargo. The pure slot-based approach is only utilised in \cite{jia_step-wise_2021, lalla-ruiz_stowage_2020} where minimisation of ballast water and hazardous cargo distances are the main objectives, respectively. 

In ocean-going RoRo shipping, where a vessel visits multiple ports, handling time is often increased by unnecessary cargo movements (shifts) \cite{ovstebo_optimization_2011,daduna_dynamic_2023,hansen_adaptive_2020}.
Shifting involves the unnecessary movement of cargo to facilitate the loading or unloading of other cargo. Both \cite{hansen_adaptive_2020} and \cite{ovstebo_optimization_2011} aim to directly minimise shifting in their stowage plans, while \cite{puisa_optimal_2021} requires that no shifting can occur.

The stability of a stowage plan is complex and accounted for in \cite{ovstebo_optimization_2011, puisa_optimal_2021, lalla-ruiz_stowage_2020, sun_semi-automated_2022}. Transverse and vertical centres of gravity are ensured to be within their limits in \cite{ovstebo_optimization_2011}. In \cite{puisa_optimal_2021}, the notion of vessel stability is also extended to include the longitudinal centre of gravity. To facilitate finding seaworthy stowage plans, the use of ballast water tanks was studied in \cite{lalla-ruiz_stowage_2020, sun_semi-automated_2022}. Recently, the impact of cargo arrival times on stowage planning has been investigated in \cite{daduna_dynamic_2023} using the grid representation.

In terms of solution methods, \cite{hansen_2d-packing_2016} and \cite{ovstebo_optimization_2011} demonstrate the necessity of heuristics when solving realistically sized instances using the lane and grid representations. A custom heuristic is introduced in \cite{ovstebo_optimization_2011}, and this approach is further enhanced with tabu search and squeaky wheel optimisation in \cite{ovstebo_routing_2011}. Both cases show that the heuristics are faster at solving realistic-sized instances than an MIP solver. Additionally, \cite{hansen_adaptive_2020} implements an Adaptive Large Neighbourhood Search heuristic with multiple destroy and repair operators to tackle the RSSP. They also illustrate that the shifting problem can be efficiently evaluated using Dijkstra's algorithm. Again, the heuristic shows it can solve the problem faster than an MIP solver. 

The literature indicates that the representation and multi-port version of the RSSP impact the choice of solution method and computational efficiency. We see a gap in the literature regarding a comparative analysis of cargo placement modelling to show its strengths and weaknesses. Each of the representations can be adapted to facilitate comparison. Furthermore, there is a need for a comprehensive model of vessel stability with ballast tanks and stress forces, similar to that for container stowage planning \cite{van_twiller_literature_2023}.

\section{Problem Formulation} \label{sec:formulation}
\begin{figure}[H]
     \centering
     \begin{subfigure}[b]{0.2\textwidth}
         \centering
         \includegraphics[width=\textwidth]{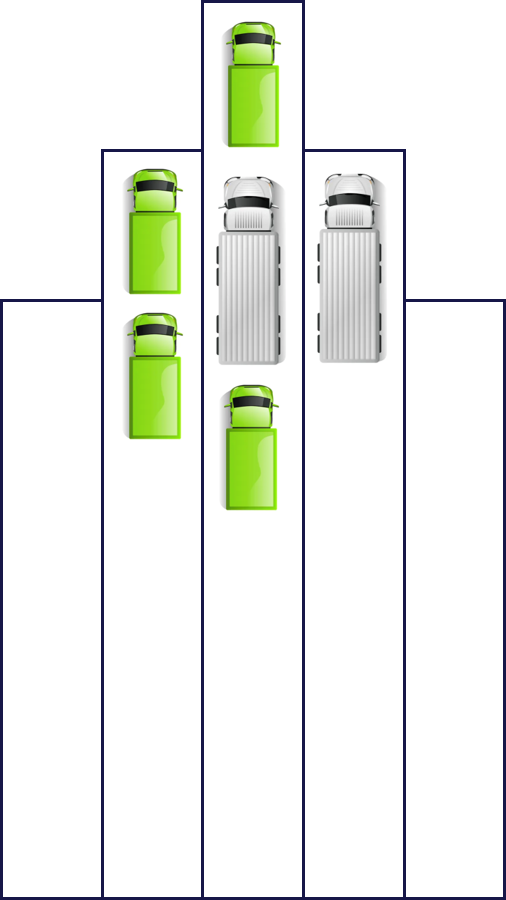}
         \caption{Lane}
         \label{fig:lane}
     \end{subfigure}
     \hfill
     \begin{subfigure}[b]{0.2\textwidth}
         \centering
         \includegraphics[width=\textwidth]{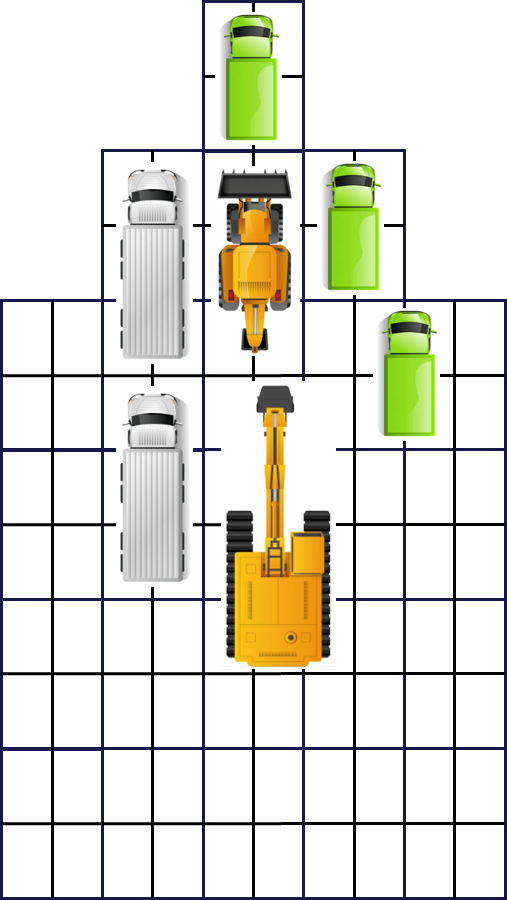}
         \caption{Grid}
         \label{fig:grid}
     \end{subfigure}
      \hfill
     \begin{subfigure}[b]{0.2\textwidth}
         \centering
         \includegraphics[width=\textwidth]{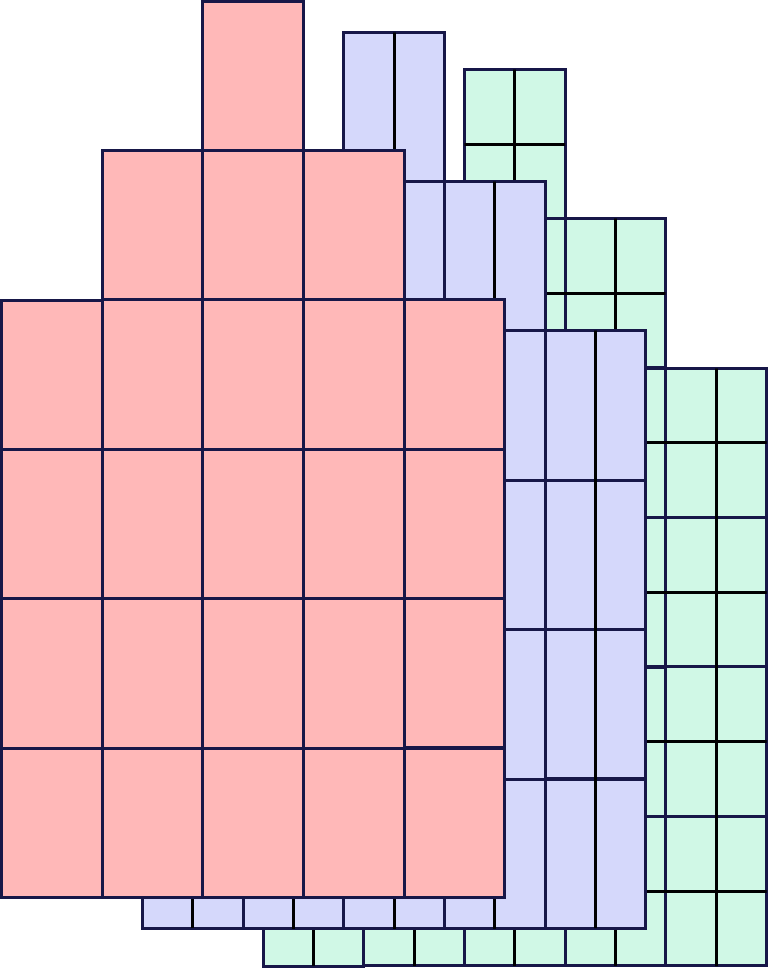}
         \caption{Hybrid-slot}
         \label{fig:hybrid}
     \end{subfigure}
        \caption{Overview of representations}
        \label{fig:representations}
\end{figure}

The literature shows different modelling approaches and objectives for the RSSP. The three representations are illustrated in \autoref{fig:representations}. We aim to quantify the strengths and weaknesses of the lane, grid, and hybrid-slot models and establish a guide for determining when each is preferable. To do that, we first need to establish the underlying principles of each representation and how they relate to the overall problem of stowage planning for RoRo vessels. A key distinguishing factor here is also what works from an Operations Research perspective and what works from an industry perspective.

\begin{figure}
    \centering
    \includegraphics[width=0.75\linewidth]{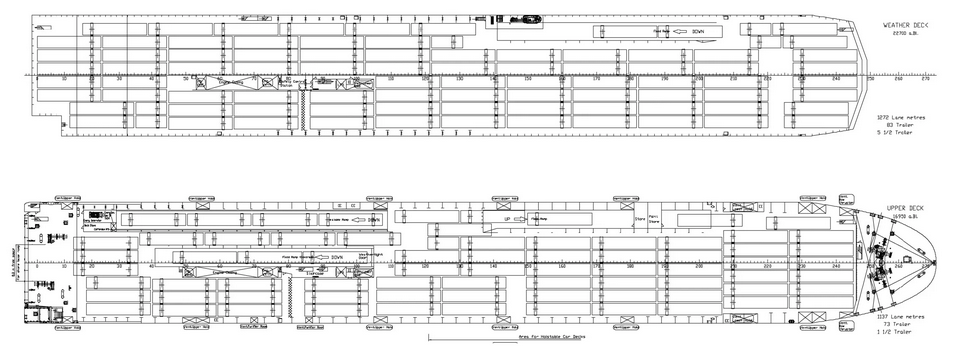}
    \caption{Overview of a RoRo vessels decks \cite{wartsila_corporation_ro-ro_2016}}
    \label{fig:deck_layout}
\end{figure}

As mentioned, the deck space of a RoRo vessel is open, and cargo can be positioned anywhere on the deck. However, the deck geometry is typically non-rectangular, featuring angled corners, central support beams, and side ramps, as shown in \autoref{fig:deck_layout}. These aspects complicate the modelling of the deck space. The lane-based representation leverages the principle of loading operations, where cargo is placed within painted lanes on the deck. This is beneficial from an industry perspective because it mimics the actual loading operation. The main problem with the lane representation lies in its mathematical modelling. Difficulties occur when cargo is wider than a single lane or when the deck is non-rectangular. Conversely, potential space is wasted when cargo is narrower than the lane width. These issues are straightforward for a human planner to handle, but pose a challenge in an optimisation model. The grid-based representation mitigates the main modelling limitations of the lane approach by allowing cargo to be placed at any grid position on the deck area. This should improve the spatial efficiency of the representation, albeit at the cost of other trade-offs. From an industry perspective, unrestricted placement may be impractical, since cargo often must be lashed or trailers positioned on  SAT trestles that must align with fixed trestle holes in the deck. On the contrary, it might be beneficial with its spatial efficiency, especially if the deck has many irregular features. Modelling-wise, it introduces the parameter of grid size, which must be chosen small enough to capture spatial detail yet large enough to remain computationally tractable. The hybrid-slot representation is the simplest in terms of mathematical modelling, as the slots are predetermined. Cargo must be placed in predefined slots, which may result in a loss of spatial efficiency, but this should make the model more computationally efficient. From an industry standpoint, it may be preferable to a purely grid-based approach because predefined slots can be aligned with deck features, and the typical way of working is more easily accommodated, albeit at the cost of requiring the slot configuration to be decided in advance. 

A further critical aspect of stowage planning is ensuring that the vessel remains seaworthy throughout the voyage. The literature presents various approaches to handling stability requirements, ranging from simplified approximations to comprehensive physical models. In this study, we evaluate full stability requirements: centres of gravity must remain within permissible limits, stress forces must not exceed structural thresholds, and ballast water adjustments are permitted to correct any imbalances. Calculating precise stability measures necessitates knowledge of the exact position and weight of each cargo item. This introduces an additional layer of complexity in the modelling of the representations, particularly for the lane-based model, where cargo positions are not directly specified, as will be discussed in subsequent sections. This formulation thus provides insight into how each representation accommodates stability constraints.

From a computational complexity standpoint, all three representations, lane, grid, and hybrid-slot, are $\mathcal{NP}$-hard, as they are variants of a multidimensional Knapsack problem \cite{ovstebo_optimization_2011, hansen_2d-packing_2016}. Our expectations and the existing literature indicate that the problem is hard to solve. We have restricted the study to Mixed-Integer Programming (MIP) solvers to ensure an equal comparison, and also to highlight the strengths and weaknesses of each model. This approach isolates the intrinsic strengths of each representation, allowing for a clear analysis of their solution quality (e.g. space utilisation) and performance across varying scales.

The literature shows, as mentioned, different modelling approaches and objectives for the RSSP. Here, the lane \cite{ovstebo_optimization_2011}, grid \cite{hansen_2d-packing_2016}, and hybrid-slot \cite{puisa_optimal_2021} approaches will be modelled consistently to enable comparison across models. This will include some changes to the previous models, particularly the lane-based model. The RSSP model presented here is set in the context of short-sea shipping, where vessels operate between two fixed ports and discharge their entire cargo at each port, called single-port planning. We accommodate heterogeneous cargoes with varying rectangular dimensions and weights. Any resulting stowage plan must be seaworthy and meet all the vessel’s stability requirements. We aim to maximise the number of cargoes carried on the voyage while minimising the ballast water required to ensure the vessel’s seaworthiness. A pseudo-revenue value is assigned to each cargo, as actual revenue data are often unavailable in practice. The model and tests presented do not distinguish between optional and required cargo; however, this could easily be introduced into the model, though it would require additional consideration in the test cases to ensure feasibility.

\section{Mathematical Formulations} \label{sec:math}
The modelling shares notation. First, we introduce how the vessel stability and stress forces are implemented in \autoref{sec:stability}. Second, we describe each of the three representations individually. The goal of the modelling is to be able to generate stowage plans that account for heterogeneous-sized cargo in their length, width and weight. All cargo is optional to load and we plan for single port voyages. The vessels hydrostatic data is assumed to have 0 trim and the hydrostatic data it self is discretised.

\subsection{Stability Requirements}
\label{sec:stability}
In this paper, we define seaworthiness as follows. The vessel's centers of gravity Vertical Center of Gravity (VCG), Transverse Center of Gravity (TCG), and Longitudinal Center of Gravity (LCG) should be within their respective limits. When the centres of gravity are within their limits, the vessel should be stable. The stress forces, shear force, and bending moment should be within the defined limits of the vessel. The limits on the stress forces ensure that the vessel's load is evenly distributed, preventing bending or flexing. The vessel's ballast tanks are modelled to correct any imbalances in weight distribution, as is done in reality. The vessel must have a minimum draft to ensure the propeller is fully immersed in the water, a requirement important  to the vessel's performance. The modelling of the seaworthiness assumes that the vessel is in calm waters, and we do not account for free surface effects as the vessels are assumed to have baffle tanks and an automated ballast management system. In addition, trim is modelled to be zero, meaning that the vessel is at even keel. The following formalises the stability requirements that are applied to the different RSSP model representations. 

\begin{table}[h!]
    \footnotesize
    \centering
        \begin{tabular}{|l|p{12cm}|}
        \hline
        \multicolumn{2}{|c|}{\textbf{SHARED SETS AND PARAMETERS}} \\ \hline
        \multicolumn{2}{|l|}{\textit{Sets}} \\ \hline
        $\mathcal{C}$ & Set of cargo units to be stowed \\
        $\mathcal{C}^n$ & Subset of cargo of type $n$ \\ 
        $\mathcal{N}^C$ & Set of cargo types \\ 
        $\mathcal{D}$ & Set of decks on the vessel \\ 
        $\mathcal{F}$ & Set of frames (longitudinal sections) of the vessel \\ 
        $\mathcal{T}$ & Set of ballast water tanks \\ 
        $\mathcal{H}$ & Set of indices in the hydrostatic table sorted by ascending weight \\ \hline
        \multicolumn{2}{|l|}{\textit{Cargo Parameters}} \\ \hline
        $W^C_c$ & Weight of cargo $c$ (tonnes) \\ 
        $L^C_c$ & Length of cargo $c$ (meters) \\ 
        $U^C_c$ & Width of cargo $c$ (meters) \\ 
        $H^C_c$ & Height of cargo $c$ (meters) \\ 
        $\pi_c$ & Pseudo-revenue for placing cargo $c$ \\ \hline
        \multicolumn{2}{|l|}{\textit{Vessel Parameters}} \\ \hline
        $W^D_d$ & Maximum weight capacity of deck $d$ (tonnes) \\ 
        $W^{\text{ship}}_f$ & Lightship weight at frame $f$ (tonnes) \\ 
        $M^L_{\text{ship}}$ & Longitudinal moment contribution from vessel and deadweights (tonne-meters) \\ 
        $M^T_{\text{ship}}$ & Transverse moment contribution from vessel and deadweights (tonne-meters) \\
        $M^V_{\text{ship}}$ & Vertical moment contribution from vessel and deadweights (tonne-meters) \\
        $L^F_f$ & Distance between consecutive frames (meters) \\
        $B_f$ & Starting position of frame $f$ (meters from vessel reference point) \\ \hline
        \multicolumn{2}{|l|}{\textit{Ballast Tank Parameters}} \\ \hline
        $T^{\max}_t$ & Maximum capacity of ballast tank $t$ (tonnes) \\
        $\text{LCG}^T_t$ & Slope of linear relationship between ballast water weight and LCG moment for tank $t$ \\ 
        $\text{TCG}^T_t$ & Slope of linear relationship between ballast water weight and TCG moment for tank $t$ \\ 
        $\text{VCG}^T_t$ & Slope of linear relationship between ballast water weight and VCG moment for tank $t$ \\
        $P^{\text{ballast}}_{f,t}$ & Fraction of ballast tank $t$ located in frame $f$ \\ \hline
        \multicolumn{2}{|l|}{\textit{Hydrostatic Parameters}} \\ \hline
        $W^{\min}_h$ & Lower bound of vessel weight for hydrostatic index $h$ (tonnes) \\
        $W^{\max}_h$ & Upper bound of vessel weight for hydrostatic index $h$ (tonnes) \\
        $\text{LCG}^{\min}_h$ & Minimum LCG limit at hydrostatic index $h$ (meters) \\
        $\text{LCG}^{\max}_h$ & Maximum LCG limit at hydrostatic index $h$ (meters) \\
        $\text{TCG}^{\min}$ & Minimum TCG limit (meters) \\ 
        $\text{TCG}^{\max}$ & Maximum TCG limit (meters) \\ 
        $\text{VCG}^{\max}_h$ & Maximum VCG limit at hydrostatic index $h$ (meters) \\ 
        $B_{h,f}$ & Cumulative buoyancy at frame $f$ for hydrostatic index $h$ (tonnes) \\ \hline
        \multicolumn{2}{|l|}{\textit{Stress Force Parameters}} \\ \hline
        $\text{SF}^{\min}_f$ & Minimum allowable shear force at frame $f$ (tonnes) \\ 
        $\text{SF}^{\max}_f$ & Maximum allowable shear force at frame $f$ (tonnes) \\ 
        $\text{BM}^{\max}_f$ & Maximum allowable bending moment at frame $f$ (tonne-meters) \\ 
        $q_f$ & Weight factor for forward vs. backward bending moment calculation at frame $f$ \\ \hline
        \end{tabular}
    \caption{Shared sets and parameters for all model formulations.}
    \label{tab:shared_notation}
\end{table}

\begin{table}[H]
    \footnotesize
    \centering
        \begin{tabular}{|l|p{11cm}|}
        \hline
        \multicolumn{2}{|c|}{\textbf{SHARED VARIABLES}} \\ \hline
        \multicolumn{2}{|l|}{\textit{Ballast Variables}} \\ \hline
        $y_t$ & Amount of ballast water in tank $t$ (tonnes) \\ \hline
        \multicolumn{2}{|l|}{\textit{Stability Variables}} \\ \hline
        $z_h$ & Binary: 1 if vessel weight falls within hydrostatic index $h$, 0 otherwise \\ 
        $m^L_{\text{cargo}}$ & Total LCG moment contribution from cargo (tonne-meters) \\ 
        $m^T_{\text{cargo}}$ & Total TCG moment contribution from cargo (tonne-meters) \\ 
        $m^V_{\text{cargo}}$ & Total VCG moment contribution from cargo (tonne-meters) \\ 
        $w^{\text{cargo}}_f$ & Total cargo weight at frame $f$ (tonnes) \\ 
        $w^{\text{ballast}}_f$ & Total ballast water weight at frame $f$ (tonnes) \\ 
        $w^{\text{total}}_f$ & Total weight (ship + cargo + ballast) at frame $f$ (tonnes) \\ 
        $b_f$ & Buoyancy at frame $f$ (tonnes) \\ 
        $sf_f$ & Shear force at frame $f$ (tonnes) \\ 
        $bm^{\text{for}}_f$ & Bending moment calculated from forward at frame $f$ (tonne-meters) \\ 
        $bm^{\text{back}}_f$ & Bending moment calculated from aft at frame $f$ (tonne-meters) \\ \hline
        \end{tabular}
    \caption{Shared variables for all model formulations.}
    \label{tab:shared_variables}
\end{table}

\begin{align}
    & \label{eq:stability_max_ballast_water}  y_t \leq T^{\max}_t &  \quad \forall t \in \mathcal{T} \\
    & \label{eq:lcg} \sum_{h \in \mathcal{H}} \text{LCG}_h^{\min} \cdot z_h \leq M_{\text{ship}}^L + m_{\text{cargo}}^L + \sum_{t \in \mathcal{T}} (y_t \cdot \text{LCG}_t^T) \leq \sum_{h \in \mathcal{H}} \text{LCG}_h^{\max} \cdot z_h \\
    & \label{eq:tcg} \sum_{h \in \mathcal{H}} \text{TCG}^{\min} \cdot z_h  \leq M_{\text{ship}}^T + m_{\text{cargo}}^T + \sum_{t \in \mathcal{T}} (y_t \cdot \text{TCG}_t^T) \leq \sum_{h \in \mathcal{H}} \text{TCG}^{\max} \cdot z_h \\
    & \label{eq:vcg} M_{\text{ship}}^V + m_{\text{cargo}}^V + \sum_{t \in \mathcal{T}} (y_t \cdot \text{VCG}^{T}_t) \leq \sum_{h \in \mathcal{H}} \text{VCG}_h^{\max} \cdot z_h \\
    & \label{eq:stability_cumulative_weight} w_f^{\text{total}} = w_{f-1}^{\text{total}} + W_f^{\text{ship}} + w_f^{\text{cargo}} + w_f^{\text{ballast}} & \quad \forall f > 0 \in \mathcal{F} \\
    & \label{eq:stability_ballast_frame} w_f^{\text{ballast}} = \sum_{t \in \mathcal{T}} P_{f,t}^{\text{ballast}} \cdot y_t & \quad \forall f \in \mathcal{F} \\
    & \label{eq:stability_dis} \sum_{h \in \mathcal{H}} (W^{\min}_h \cdot z_h) \leq w^{\text{total}}_{|\mathcal{F}|} \leq \sum_{h \in \mathcal{H}} (W^{\max}_h \cdot z_h)\\
    & \label{eq:stability_dis_one} \sum_{h \in \mathcal{H}} z_h \leq 1 \\
    & \label{eq:stability_buoy} \frac{\sum_{h \in \mathcal{H}} B_{h,f} \cdot z_h}{2}  = b_f & \quad \forall f \in \mathcal{F} \\
    & \label{eq:stability_shear} w_f^{\text{total}} - b_f = sf_f & \quad \forall f \in \mathcal{F} \\
    & \label{eq:stability_shear_limit} \text{SF}_f^{\min} \leq sf_f \leq \text{SF}_f^{\max} & \quad \forall f \in \mathcal{F} \\
    & \label{eq:stability_bending_for} bm^{\text{for}}_f = \frac{sf_{f-1} + sf_f}{2} \cdot L^F_f + bm^{\text{for}}_{f-1} & \quad \forall f \in \mathcal{F} \\
    & \label{eq:stability_bending_back} bm^{\text{back}}_f = \frac{sf_{f+1} + sf_f}{2} \cdot L^F_f + bm^{\text{back}}_{f+1} & \quad \quad \forall f \in \mathcal{F} \\
    & \label{eq:stability_bending_limit} q_f \cdot bm^{\text{for}}_f + (1-q_f) \cdot bm^{\text{back}}_f \leq \text{BM}_f^{\max} & \quad \forall f \in \mathcal{F} \\
    & \label{eq:stability_variable_domain_continuous} y_t \geq 0 \in \mathbb{R} & \quad \forall t \in \mathcal{T} \\
    & \label{eq:stability_variable_domain_binary}  z_h \in \{0,1\} & \quad \forall h \in \mathcal{H} \\
    & \label{eq:stability_variable_domain_cg}  m^{L}_{\text{cargo}},  m^{T}_{\text{cargo}},  m^{V}_{\text{cargo}} \in \mathbb{R}  \\ 
    & \label{eq:stability_variable_domain_w}  w_f, w^{\text{cargo}}_f, w^{\text{ballast}}_f \in \mathbb{R} & \quad \forall f\in \mathcal{F} \\
    & \label{eq:stability_variable_domain_w}  b_f, sf_f, bm_f\in \mathbb{R} & \quad \forall f\in \mathcal{F}\\ \nonumber
\end{align}

\begin{figure}
    \centering
    \begin{subfigure}[t]{0.45\textwidth}
    \begin{tikzpicture}
        \begin{axis}[
            width=6cm,  
            height=6cm,  
            xlabel=Weight,
            ylabel=VCG,
            grid=both,
            legend pos = north west 
        ]
        
        \draw[black, dashed, thick] (axis cs:\pgfkeysvalueof{/pgfplots/xmin},0) -- 
                                     (axis cs:\pgfkeysvalueof{/pgfplots/xmax},0);
                
        \addplot[red, thick] table[x=Weight, y=VCG] {ballast.dat};
        \addlegendentry{VCG}
        
        \end{axis}
    \end{tikzpicture}
    \caption{The weight and VCG relationship of a specific ballast tank; highlighting a big change between the minimum and maximum, but also show that the relationship is not linear. }
    \label{fig:ballast_vcg}
    \end{subfigure}
    \hfill
    \begin{subfigure}[t]{0.45\textwidth}
    \begin{tikzpicture}
        \begin{axis}[
            width=6cm,  
            height=6cm,  
            xlabel=Weight,
            ylabel=VCG $\cdot$ Weight,
            grid=both,
            legend pos = north west 
        ]
        
        \draw[black, dashed, thick] (axis cs:\pgfkeysvalueof{/pgfplots/xmin},0) -- 
                                     (axis cs:\pgfkeysvalueof{/pgfplots/xmax},0);
                
        \addplot[red, thick] table[x=Weight, y=M_VCG] {ballast.dat};
        \addlegendentry{True}

    \addplot [
        domain=0:1000, 
        samples=100, 
        color=blue,
        dotted, thick
        ]
        {x*5.4};
    \addlegendentry{Estimated}
            
        \end{axis}
    \end{tikzpicture}
    \caption{Linear approximation of the VCG moment for a ballast tank, showing the true relationship versus the estimated linear model used in the optimization.}
    \label{fig:ballast_vcg_estimated}
    \end{subfigure}
    \caption{The VGG of a ballast tank depends on its filling level. The figure show the relationship and how it is estimated linearly.}
\end{figure}
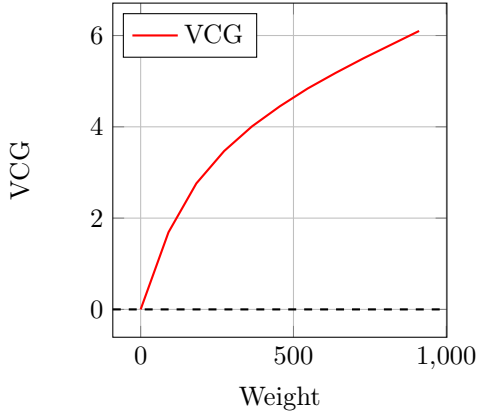
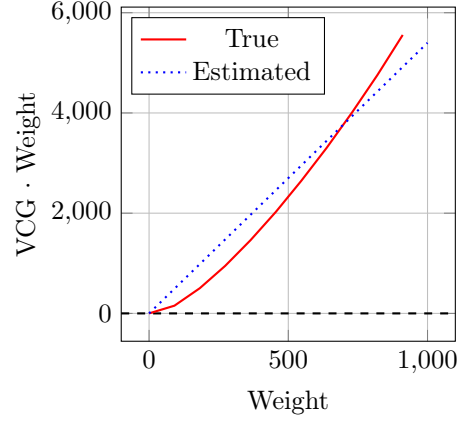

The stability calculations represented in equations \eqref{eq:stability_max_ballast_water} to \eqref{eq:stability_variable_domain_w} are used across all three model representations shown in this paper. However, cargo placement varies across models; hence, cargo-related stability variables are specific to each model. All cargo related stability variables (with the subscript \textit{cargo}) are detailed in \autoref{sec:formulation} for each of the three models. The first constraints \eqref{eq:stability_max_ballast_water} ensure that the amount of ballast water each tank $t\in T$ is within the capacity $T^{\max}_t$ of the tank $t$. Constraints \eqref{eq:lcg}-\eqref{eq:vcg} ensure that the LCG, TCG and VCG of the vessel are within their limits. The LCG of the vessel in constraint \eqref{eq:lcg} is the sum of the contribution to the LCG from the vessel itself and deadweights such as the fuel tanks represented by $M^L_{\text{ship}}$, the cargo represented by $m^L_{\text{cargo}}$ and the ballast tanks represented by $y_t \cdot LCG^T_t$. The LCG limits are found from the Longitudinal Center of Buoyancy (LCB). Since we require the trim of the vessel to be zero, we can use the principle that a vessel achieves stable equilibrium (zero trim) when the LCG is equal to the LCB. The LCB is dynamic and depends on the vessel's weight. To find the value of LCB, we use the selector variable $z_h$ to find an upper $\text{LM}_h^{\max}$ and lower $\text{LM}_h^{\min}$ limit on the LCB. $z_h$ selects the hydrostatic index that contains the current total weight. How the weight interval $z_h$ is found is explained in the description of constraints \eqref{eq:stability_dis} and \eqref{eq:stability_dis_one}. Constraint \eqref{eq:tcg} limits the TCG of the vessel. The TCG of the vessel is not dynamic and should be zero to prevent undesired heel. The bounds $\text{TM}^{\min}$ and $\text{TM}^{\max}$ are therefore a lower and upper bound of zero that are within the limits of the vessel. Constraint \eqref{eq:vcg} limits the VCG og the vessel. The VCG is dynamic and depends on the Metacentric Height ($GM$) and maximum metacentre ($KMT$). The $GM$ of a vessel is given by $GM=KMT-VCG$. The hydrostatic table provides a minimum $GM_h$ and maximum $KMT_h$ for each index $h$ and is used to calculate $V^{\max}_h = KMT_h - GM_h$. The VCG of a ballast tank depends on its filling level, shown in \autoref{fig:ballast_vcg}. Instead of discretising the ballast tank filling levels as in \cite{lalla-ruiz_stowage_2020}, the VCG contribution of the ballast tanks is estimated here using a linear approximation. Specifically, $\text{VCG}^{T}_t$ represents the slope of the linear relationship between ballast water weight and the VCG moment as seen in \autoref{fig:ballast_vcg_estimated}, allowing us to calculate the VCG contribution. The contributions for the LCG and the TCG exhibit the same problem to a lesser extent when the tanks are non-cubic and modelled similarly.

The cumulative weight $w^{total}_f$ of each frame $f \in \mathcal{F}$ of the vessel is calculated in \eqref{eq:stability_cumulative_weight}. Here, $W^{\text{ship}}_f$ is the lightship weight, $w^{\text{cargo}}_f$ is the cargo weight and $W^{\text{ballast}}_f$ is the ballast weight in each frame respectively. The total ballast water weight at each frame, $w^{ballast}_f$, is given by constraint \eqref{eq:stability_ballast_frame} where the parameter $P^{tank}_{f,t}$ is the fraction of the tank $t$ in frame $f$ (this is necessary as tanks and frames do not necessarily align). As previously explained, the stability limits are found in the vessel's hydrostatic data table. Assuming the vessel needs to sail at even-keel (at a trim of zero), it is possible to obtain those limits by simple interpolation \cite{pacino_accurate_2012}. The interpolation is based on the vessel's total weight. We introduce the set $\mathcal{H}$ where each $h$ is an index in the hydrostatic table sorted by ascending weight. Let $W^{\min}_h$ and $W^{\max}_h$ be the lower and upper bounds of the total vessel weight at each hydrostatic index. Constraint \eqref{eq:stability_dis} force the total weight of the vessel $w^{\text{total}}_{|F|}$, to be within an index using the binary indicator variable $z_h$. We then impose that only one index can be active in constraint \eqref{eq:stability_dis_one}. Further, the set $\mathcal{H}$ is limited to only include indices that ensure solutions with the minimum draft required by the propeller. This ensures that the vessel's total weight is between exactly two adjacent entries represented by $h\in \mathcal{H}$ in the hydrostatic table. Constraint \eqref{eq:stability_buoy} utilises the hydrostatic indices to interpolate the cumulative buoyancy $B_{h,f}$ at each frame. The interpolation is done for each frame because it is needed for the stress forces calculation. In constraint \eqref{eq:stability_shear}, the stress forces of each frame $sf_f$ are calculated as the difference between the accumulated weight and the accumulated buoyancy at the frame. The shear force is then constrained between a lower $SF^{\min}_f$ and upper $SF^{\max}_f$ limit for each frame in \eqref{eq:stability_shear_limit}. Finally, the accumulated bending moment $bm_f$, at each frame, is calculated as the area under the curve defined by the shear force $sf_f$. This is done in constraints \eqref{eq:stability_bending_for} and \eqref{eq:stability_bending_back}, where the moment is calculated from each side of the vessel. Here, $L^F_f$ is the distance between two consecutive frames. The bending moment is then constrained to be within its limits ($BM^{max}_f$) with constraint \eqref{eq:stability_bending_limit} where $q_f$ is the fraction of how much the back and forward bending are contributing to the bending moment. The domains of the variables are defined in \eqref{eq:stability_variable_domain_continuous} to \eqref{eq:stability_variable_domain_w}. The stability constraints will remain the same in all three models described in the following sections. 

\subsubsection{Hybrid formulation} \label{sec:formulation_hybrid}
The hybrid-slot representation divides the deck layout into slots for each cargo type that needs to be planned. To ensure that cargo does not overlap, the slots that utilise the same area are identified. The set of overlapping slots can be used to constrain which slots cannot be used simultaneously. The slots and overlaps are provided from vessel data or pre-computed and are therefore not part of the optimisation. The hybrid-slot representation was introduced in \cite{puisa_optimal_2021} and is here extended to include ballast-water and stress-force constraints, as discussed in \autoref{sec:stability}. The mathematical formulation is given below, and the notation in \autoref{tab:hybrid_notation}.

\begin{table}[H]
    \footnotesize
    \centering
        \begin{tabular}{|l|p{7cm}|}
        \hline
        \multicolumn{2}{|c|}{\textbf{HYBRID-SLOT MODEL}} \\ \hline
        \multicolumn{2}{|l|}{\textit{Sets}} \\ \hline
        $\mathcal{S}$ & Set of all pre-computed slots on the vessel \\ 
        $\mathcal{S}^n$ & Subset of slots suitable for cargo type $n$ \\ 
        $\mathcal{S}^D_d$ & Subset of slots on deck $d$ \\ 
        $\mathcal{S}^F_f$ & Subset of slots intersecting frame $f$ \\ 
        $\mathcal{O}$ & Set of overlapping slot pairs $(i,j)$ \\ \hline
        \multicolumn{2}{|l|}{\textit{Parameters}} \\ \hline
        $\text{LCG}^S_s$ & Longitudinal center of gravity of slot $s$ (meters) \\ 
        $\text{TCG}^S_s$ & Transverse center of gravity of slot $s$ (meters) \\ 
        $\text{VCG}^S_s$ & Vertical center of gravity of slot $s$ (meters) \\ 
        $P^S_{f,s}$ & Fraction of slot $s$ located in frame $f$ \\ \hline
        \multicolumn{2}{|l|}{\textit{Variables}} \\ \hline
        $x_{c,s}$ & Binary: 1 if cargo $c$ is assigned to slot $s$, 0 otherwise \\ \hline
        \end{tabular}
    \caption{Hybrid-slot model specific notation.}
    \label{tab:hybrid_notation}
\end{table}


\begin{align}
    \label{eq:hybrid_objective}
    \text{Min Z = } \quad & \sum_{t \in \mathcal{T}} y_t - \sum_{c \in \mathcal{C}}\sum_{s \in \mathcal{S}}  \pi_c \cdot x_{c,s} \\ 
    \text{s.t.} & \nonumber \\
    & \label{eq:hybrid_slot_limit} \sum_{c \in \mathcal{C}} x_{c, s} \leq 1 & \quad  \forall s \in \mathcal{S} \\
    & \label{eq:hybrid_slot_to_cargo} \sum_{s \in \mathcal{S}^n} x_{c, s} \leq 1 \quad &  \forall c \in \mathcal{C}^n, n \in \mathcal{N^C} \\
    & \label{eq:hybrid_slot_overlap} \sum_{c \in \mathcal{C}} x_{c,i} + \sum_{c \in \mathcal{C}} x_{c,j}  \leq 1 \quad & \forall (i, j) \in \mathcal{O} \\
    & \label{eq:hybrid_deck_weight} \sum_{c \in \mathcal{C}} \sum_{s\in\mathcal{S}^D_d} W^C_c \cdot x_{c,s} \leq W^{D}_d \quad & \forall d \in \mathcal{D} \\
    & \label{eq:hybrid_slot_cargo_frame} \sum_{c \in \mathcal{C}} \sum_{s \in \mathcal{S}^F_f} P^{S}_{f,s} \cdot W^C_c \cdot x_{c,s} = w^{cargo}_f \quad &  \forall f \in \mathcal{F} \\
    & \label{eq:hybrid_lcg_cargo} \sum_{c \in \mathcal{C}} W^C_c \cdot x_{c,s} \cdot LCG^{S}_s = m^{L}_{cargo} \\
    & \label{eq:hybrid_tcg_cargo} \sum_{c \in \mathcal{C}} W^C_c \cdot x_{c,s} \cdot TCG^{S}_s = m^{T}_{cargo} \\
    & \label{eq:hybrid_vcg_cargo} \sum_{c \in \mathcal{C}} W^C_c \cdot x_{c,s} \cdot VCG^{S}_s = m^{V}_{cargo} \\
    & \label{eq:hybrid_variable domain_a} x_{c,s} \in \{0,1\} \quad & \forall c \in C, \forall s \in S \\
    & \label{eq:hybrid_variable domain_cg}  m^{L}_{cargo},  m^{T}_{cargo},  m^{V}_{cargo},y_t \in \mathbb{R} \\
    & \text{Stability is subject to equations \eqref{eq:stability_max_ballast_water} - \eqref{eq:stability_variable_domain_w}} \nonumber 
\end{align}

The objective function \eqref{eq:hybrid_objective} minimises the water in all ballast tanks while maximising the cargo placed in the stowage plan. $y_t$ is the amount of ballast water in each ballast tank $t\in\mathcal{T}$. The binary variable $x_{c,s}$ is the assignment of cargo $c\in \mathcal{C}$ to slot $s\in \mathcal{S}$, where the set $\mathcal{S}$ contains all the pre-computed slots on the vessel. The parameter $\pi_c$ is the pseudo-revenue for placing cargo. We assume that bringing more cargo is always preferable to using less ballast water; thus, $\pi_c$ is equal to the total amount of ballast water possible to ensure that we weigh bringing cargo higher than the ballast water usage. Constraints \eqref{eq:hybrid_slot_limit} ensure that each slot will accommodate at most one cargo. Slots can only accommodate cargo of a specific type, which is modelled by constraints  \eqref{eq:hybrid_slot_to_cargo}. Here $n\in \mathcal{N}^c$ defines the cargo type, and the sets  $\mathcal{C}^n$ and $\mathcal{S}^n$ are the subsets of cargo and slots available for cargo type $n$, respectively. Overlapping slots are avoided in constraints \eqref{eq:hybrid_slot_overlap} by ensuring that for all (unordered) overlapping slot pairs $(s_i,s_j)\in \mathcal{O}$ at most one slot can be selected. The deck weight limits are constrained in \eqref{eq:hybrid_deck_weight} by the summation of the weight of all slots on the deck $W^D_d$ for all decks $d \in \mathcal{D}$. Constraints \eqref{eq:hybrid_slot_cargo_frame} determine the cargo weight in each frame of the vessel, denoted $w^{cargo}_f$. The set $\mathcal{S}^F_f$ is all the slots in a frame $f\in \mathcal{F}$ and the parameter $P^{S}_{f,s}$ is the fraction of slot $s$ in the frame $f$. Finally, constraints \eqref{eq:hybrid_lcg_cargo}-\eqref{eq:hybrid_vcg_cargo} compute the centre of gravity needed by the stability constraints from \autoref{sec:stability}. Each slot $s$ has its own pre-calculated $LCG_s$, $TCG_s$, and $VCG_s$, centred in the middle of the slot. The cargo weight in each slot $s$ is multiplied by the corresponding $LCG_s$, $TCG_s$, and $VCG_s$ to calculate the total cargo contribution. The variables $m^{L}_{cargo}$, $m^{T}_{cargo}$ and $m^{V}_{cargo}$ are the total contributions of the cargo in the longitudinal, transverse, and vertical direction, respectively.

\subsubsection{Grid formulation} \label{sec:formulation_grid}
The grid representation allows for arbitrary assignment of cargo anywhere on the deck. The deck layout is split into a grid-based coordinate system for placing cargo, similar to the modelling of 2D packing models. Cargo is placed based on the lower-left corner of its bounding box. The complete set of stability constraints and the modelling of ballast water have been added to the original model \cite{hansen_2d-packing_2016}. 
As before, a summary of the notation can be found in \autoref{tab:grid_notation}. The mathematical formulation is given below.

\begin{table}[H]
    \footnotesize
    \centering
        \begin{tabular}{|l|p{10cm}|}
        \hline
        \multicolumn{2}{|c|}{\textbf{GRID MODEL}} \\ \hline
        \multicolumn{2}{|l|}{\textit{Sets}} \\ \hline
        $\mathcal{I}$ & Set of grid rows (longitudinal positions) across all decks \\ 
        $\mathcal{J}$ & Set of grid columns (transverse positions) across all decks \\ 
        $\mathcal{I}^c_i$ & Set of grid rows that cargo $c$ occupies when placed at row $i$ \\ 
        $\mathcal{J}^c_j$ & Set of grid columns that cargo $c$ occupies when placed at column $j$ \\ 
        $\mathcal{I}^d$ & Subset of grid rows on deck $d$ \\ 
        $\mathcal{J}^d$ & Subset of grid columns on deck $d$ \\ 
        $\mathcal{P}^f$ & Set of grid positions $(i,j)$ within frame $f$ \\ \hline
        \multicolumn{2}{|l|}{\textit{Parameters}} \\ \hline
        $\text{LCG}^G_i$ & Longitudinal center of gravity of grid cell in row $i$ (meters) \\
        $\text{TCG}^G_j$ & Transverse center of gravity of grid cell in column $j$ (meters) \\ 
        $\text{VCG}^G_i$ & Vertical center of gravity of grid cell in row $i$ (meters) \\
        $P^G_{f,i,c}$ & Fraction of cargo $c$ in frame $f$ when placed at grid row $i$ \\ \hline
        \multicolumn{2}{|l|}{\textit{Variables}} \\ \hline
        $x_{c,i,j}$ & Binary: 1 if lower-left corner of cargo $c$ is placed at grid cell $(i,j)$, 0 otherwise \\ \hline
        \end{tabular}
    \caption{Grid model specific notation.}
    \label{tab:grid_notation}
\end{table}


\begin{align}
    \label{eq:grid_objective}
    \text{Min Z = } \quad & \sum_{t \in \mathcal{T}} y_t - \sum_{c \in \mathcal{C}} \sum_{i \in \mathcal{I}} \sum_{j \in \mathcal{J}} \pi_c \cdot x_{c,i,j} \\ 
    \text{s.t.} & \nonumber \\
    & \label{eq:grid_no_overlap} \sum_{c \in \mathcal{C}} \sum_{i' \in \mathcal{I}^c_i} \sum_{j' \in \mathcal{J}^c_j} x_{c,i',j'} \leq 1 \quad & \forall i \in \mathcal{I}, j \in \mathcal{J} \\
    & \label{eq:grid_corner} \sum_{i' \in I^c_i} \sum_{j' \in J^c_j} x_{c,i',j'} \geq x_{c,i,j}  \quad & \forall c \in \mathcal{C}, i \in \mathcal{I}, j \in \mathcal{J} \\
    & \label{eq:grid_cumulative_weight} \sum_{(i,j) \in \mathcal{P}^f} \sum_{c \in \mathcal{C}} x_{c,i,j} \cdot W^C_c \cdot P^G_{f,i,c} = w_f^{\text{cargo}} \quad & \forall f \in \mathcal{F} \\
    & \label{eq:grid_deck_weight} \sum_{i \in \mathcal{I}^d} \sum_{j \in \mathcal{J}^d}  \sum_{c \in \mathcal{C}} x_{c,i,j} \cdot W^C_c \leq W^{D}_d \quad & \forall d \in \mathcal{D} \\
    & \label{eq:grid_lcg} \sum_{c \in \mathcal{C}} \sum_{i \in \mathcal{I}} \sum_{j \in \mathcal{J}} x_{c,i,j} \cdot W^C_c \cdot (LCG^{G}_i + \frac{L^C_c}{2}) = m^{L}_{\text{cargo}} \\
    & \label{eq:grid_tcg} \sum_{c \in \mathcal{C}} \sum_{i \in \mathcal{I}} \sum_{j \in \mathcal{J}} x_{c,i,j} \cdot W^C_c \cdot (TCG^{G}_j + \frac{U^C_c}{2}) = m^{T}_{\text{cargo}} \\
    & \label{eq:grid_vcg} \sum_{c \in \mathcal{C}} \sum_{i \in \mathcal{I}} \sum_{j \in \mathcal{J}} x_{c,i,j} \cdot W^C_c \cdot (VCG^{G}_{i} + \frac{H^C_c}{2}) = m^{V}_{\text{cargo}} \\
    & \label{eq:grid_variable_domain_a} x_{c,i,j} \in \{0,1\} \quad & \forall c \in \mathcal{C}, i \in \mathcal{I}, j \in \mathcal{J} \\
    & \label{eq:grid_variable_domain_cg} m^{L}_{\text{cargo}}, m^{T}_{\text{cargo}}, m^{V}_{\text{cargo}},y_t \in \mathbb{R} \\
    & \text{Stability is subject to equations \eqref{eq:stability_max_ballast_water} - \eqref{eq:stability_variable_domain_w}} \nonumber
\end{align}
    
The grid and hybrid-slot models are similar. Therefore, only the differences will be explained. The main variable $x_{c, i, j}$ is binary and tracks if the lower left corner of cargo $c$ is assigned to a grid cell $(i,j)$. The sets $\mathcal{I}$ and $\mathcal{J}$ are the rows and columns of the grid across all decks. The objective function \eqref{eq:grid_objective} minimises ballast water usage and maximises cargo placement, as in the hybrid model. Constraints \eqref{eq:grid_no_overlap} ensure that cargo cannot occupy the same area. The sets $\mathcal{I}^c_i$ and $\mathcal{J}^c_j$ contain the grid cells that cargo $c$ will occupy if placed at $(i,j)$. Constraints \eqref{eq:grid_corner} ensure that cargo $c$ is placed exactly in the lower left corner $(i,j)$ and not in any other grid cell it occupies. Constraints \eqref{eq:grid_cumulative_weight} determine the cargo weight in each frame of the vessel, where $\mathcal{P}^f$ contains all grid positions $(i,j)$ within the frame $f$. The parameter $P^G_{f,i,c}$ is the fraction of the cargo $c$ in frame $f$ if placed at grid $i$. Constraints \eqref{eq:grid_deck_weight} ensure the weight limit of each deck, where $\mathcal{I}^d$ and $\mathcal{J}^d$ are the grid cells on deck $d$. Constraint \eqref{eq:grid_lcg} determines the LCG moment of the cargo. Here, the LCG of the cargo is found from the LCG of the lower left grid cell on which the cargo is placed, where half of the length of the cargo $\frac{L^C_c}{2}$ is added. This ensures that the LCG of the cargo is located at its centre. The TCG and VCG of the cargo are determined in a similar fashion in \eqref{eq:grid_tcg} and \eqref{eq:grid_vcg}.

\subsubsection{Lane formulation} \label{sec:formulation_lane}
\begin{figure}[!h]

    \begin{subfigure}{.5\textwidth} 
        \centering
        \includegraphics[width=.9\linewidth]{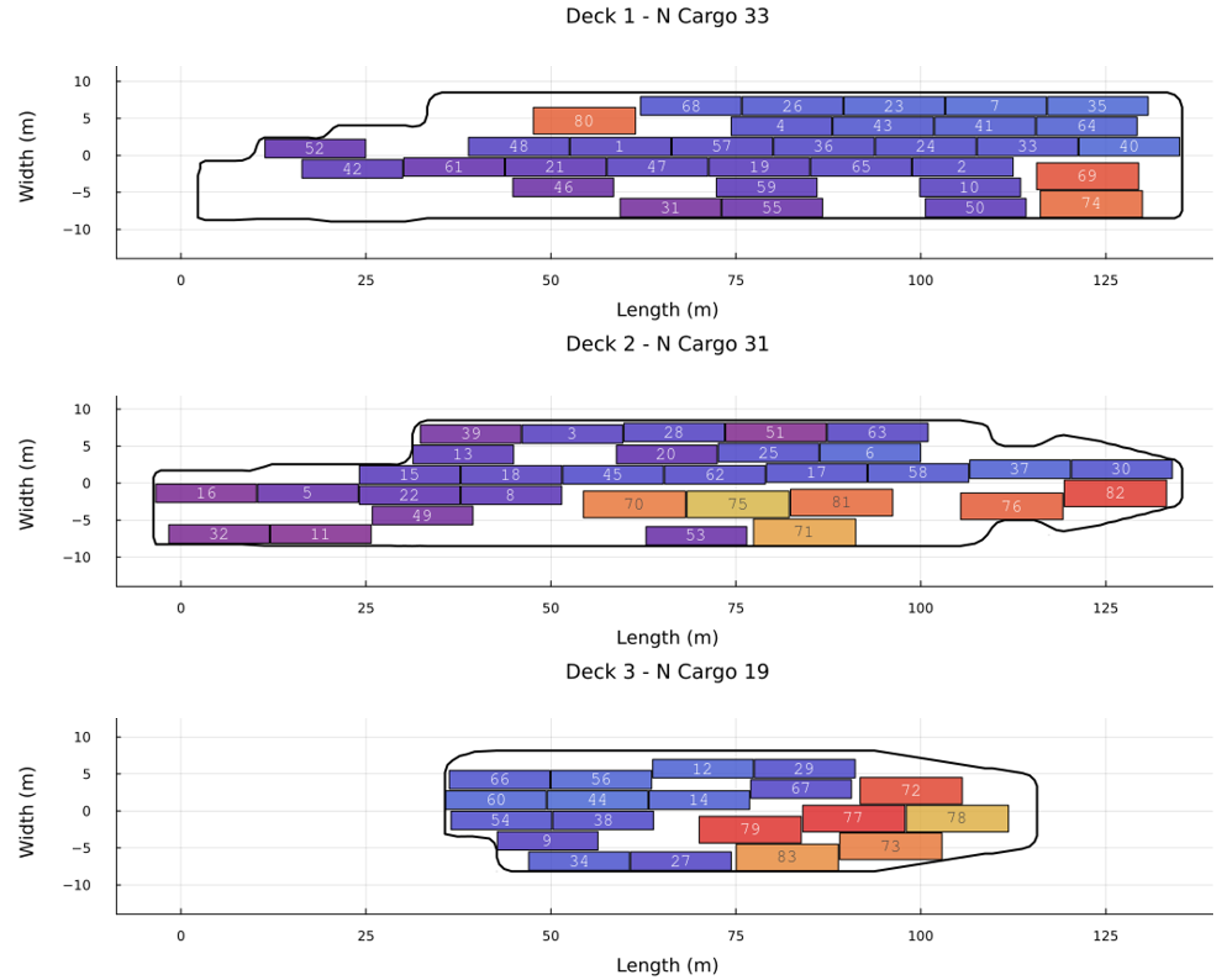}  
        \caption{Hybrid-Slot Solution}
        \label{fig:hybrid-solution}
    \end{subfigure}
    \begin{subfigure}{.5\textwidth}
        \centering
        \includegraphics[width=.9\linewidth]{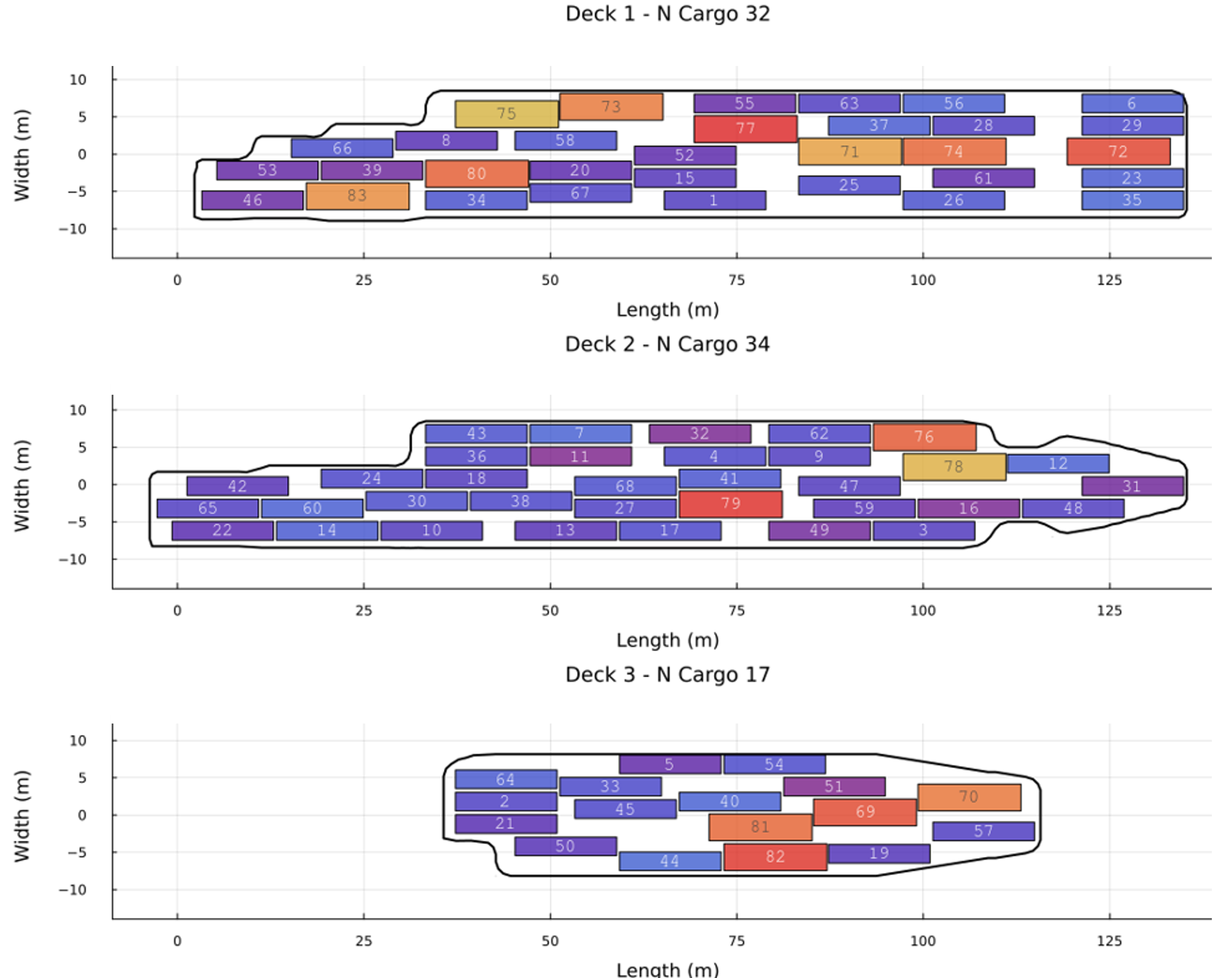}  
        \caption{Grid Solution}
        \label{fig:grid-solution}
    \end{subfigure}

    \label{fig:fig}
    \begin{subfigure}{\textwidth}
        \centering
        \includegraphics[width=.45\linewidth]{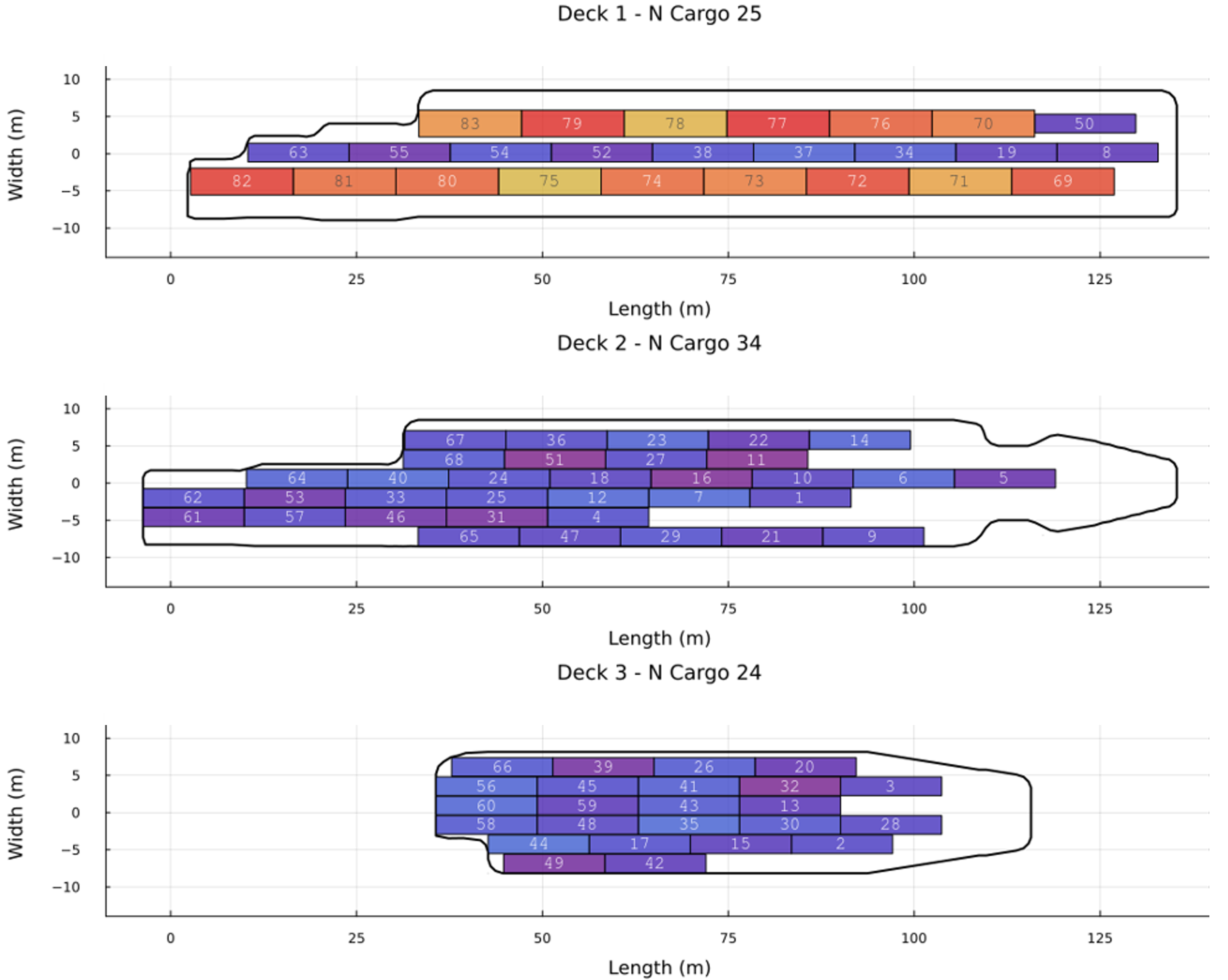}   
        \caption{Lane Solution}
        \label{fig:lane-solution}
    \end{subfigure}
    \caption{Overview of an optimal solution for all model on the small vessel. The instance is trailer and secubox, 80\% cargo and heavy cargo. All models use the same amount of ballast water 1904 Tons and places all cargo.}
    \label{fig:solutions}
\end{figure}

The lane representation is inspired by the painted lanes on the decks of RoRo vessels. Cargo is assigned to lanes, and the model keeps track of the sequence of cargo placement within each lane. To allow the modelling of complex deck layouts as seen in \autoref{fig:solutions}, we allow a predefined set of lanes to be merged to accommodate cargo wider than a standard lane width. Lanes can be merged if they are adjacent and have the same length within a predetermined tolerance. The lane-based formulation builds on the model from \cite{ovstebo_optimization_2011}. The original formulation allowed for variable lane width but also assumed that all lanes have the same length. The use of merged lanes is our modelling approach to relax the equal-lane assumption of the original model. Further modifications were made to facilitate accurate modelling of stability constraints. It should be noted that in \cite{ovstebo_optimization_2011}, the objective was to minimise cargo shifting across multiple port visits, which is very different from the problem of maximising the amount of shipped cargo while minimising ballast water for single-port stowage planning. The mathematical formulation is given below, and its notation is summarised in \autoref{tab:lane_notation}.

\begin{table}[H]
    \footnotesize
    \centering
        \begin{tabular}{|l|p{11cm}|}
        \hline
        \multicolumn{2}{|c|}{\textbf{LANE MODEL}} \\ \hline
        \multicolumn{2}{|l|}{\textit{Sets}} \\ \hline
        $\mathcal{L}$ & Set of all lanes (including single and merged lanes) across all decks \\
        $\mathcal{M}$ & Set of merged lane configurations $(p, i, j)$ where $p$ is merged lane, $i,j$ are parent lanes \\ \hline
        \multicolumn{2}{|l|}{\textit{Parameters}} \\ \hline
        $L^{\text{lane}}_l$ & Length of lane $l$ (meters) \\ 
        $U^{\text{lane}}_l$ & Width of lane $l$ (meters) \\ 
        $\text{LCG}^{\text{start}}_l$ & LCG of the starting position of lane $l$ (meters) \\ 
        $\text{TCG}_l$ & Transverse center of gravity of lane $l$ (meters) \\ 
        $\text{VCG}_l$ & Vertical center of gravity (deck level) of lane $l$ (meters) \\ 
        $M^{\max}_l$ & Maximum number of cargo units that can be assigned to lane $l$ \\ 
        $M$ & Large constant (big-M) \\ \hline
        \multicolumn{2}{|l|}{\textit{Variables}} \\ \hline
        $x_{l,c}$ & Binary: 1 if cargo $c$ is assigned to lane $l$, 0 otherwise \\ 
        $z_{l,c,c'}$ & Binary: 1 if cargo $c$ precedes cargo $c'$ in lane $l$, 0 otherwise \\ 
        $\mu_p$ & Binary: 1 if merged lane $p$ is activated, 0 otherwise \\
        $\phi_{c,f}$ & Binary: 1 if cargo $c$ is present in frame $f$, 0 otherwise \\
        $s_{l,c}$ & Starting position of cargo $c$ in lane $l$ (meters from lane start) \\ 
        $r_{c,f}$ & Starting position of cargo $c$ within frame $f$ (meters) \\ 
        $\rho_{c,f}$ & Fraction of cargo $c$ located in frame $f$ \\ 
        $\text{lcg}_c$ & Longitudinal center of gravity of cargo $c$ (meters) \\ \hline
        \end{tabular}
    \caption{Lane model specific notation.}
    \label{tab:lane_notation}
\end{table}

\begin{align}
    \label{eq:lane_objective} \text{Min Z = } \quad & \sum_{t \in \mathcal{T}} y_t - \sum_{c \in \mathcal{C}} \sum_{l \in \mathcal{L}} \pi_c \cdot x_{l,c} \\
    \text{s.t.} & \nonumber \\
    & \label{eq:lane_precedence_required} 
    z_{l,c,c'} + z_{l,c',c} \geq x_{l,c} + x_{l,c'} - 1 \quad & \forall l \in \mathcal{L}, (c,c') \in \mathcal{C}, c < c' \\
    & \label{eq:lane_precedence_exclusive} 
    z_{l,c,c'} + z_{l,c',c} \leq 1 \quad & \forall l \in \mathcal{L}, (c,c') \in \mathcal{C}, c < c' \\
    & \label{eq:lane_precedence_bound_out} 
    \sum_{c' \in \mathcal{C}} z_{l,c,c'} \leq M^{\max}_l \cdot x_{l,c} \quad & \forall l \in \mathcal{L}, c \in \mathcal{C} \\
    & \label{eq:lane_precedence_bound_in} 
    \sum_{c' \in \mathcal{C}} z_{l,c',c} \leq M^{\max}_l \cdot x_{l,c} \quad & \forall l \in \mathcal{L}, c \in \mathcal{C} \\
    & \label{eq:lane_no_overlap_forward} 
    s_{l,c'} \geq s_{l,c} + L^{\text{cargo}}_c - L^{\text{lane}}_l \cdot (1 - z_{l,c,c'}) \quad & \forall l \in \mathcal{L}, (c,c')\in \mathcal{C}, c \neq c' \\
    & \label{eq:lane_position_bound} 
    s_{l,c} \leq (L^{\text{lane}}_l - L^{\text{cargo}}_c) \cdot x_{l,c} \quad & \forall l \in \mathcal{L}, c \in \mathcal{C} \\
    & \label{eq:lane_merge_enable_new} 
    \sum_{c \in \mathcal{C}} x_{p,c} \leq M \cdot \mu_p \quad & \forall p \in \mathcal{M} \\
    & \label{eq:lane_merge_disable_original} 
    x_{i,c} + x_{j,c} \leq 2 \cdot (1 - \mu_p) \quad & \forall (p, i, j) \in \mathcal{M}, c \in \mathcal{C} \\
    & \label{eq:lane_width_constraint} 
    U^{\text{cargo}}_c \cdot x_{l,c} \leq U^{\text{lane}}_l \quad & \forall l \in \mathcal{L}, c \in \mathcal{C} \\
    & \label{eq:lane_cargo_lcg} \sum_{l \in \mathcal{L}} \left(s_{l,c} + x_{l,c} \cdot \left(L^{\text{cargo}}_c/2 + \text{LCG}_l^{\text{start}}\right)\right) = \text{lcg}_c \quad & \forall c \in \mathcal{C} \\
    & \label{eq:lane_total_lcg_moment} \sum_{c \in \mathcal{C}} \text{lcg}_c \cdot W_c = m^L_{\text{cargo}} \\
    & \label{eq:lane_total_tcg_moment} \sum_{l \in \mathcal{L}} \sum_{c \in \mathcal{C}} \text{TCG}_l \cdot W_c \cdot x_{l,c} = m^T_{\text{cargo}}  \\
    & \label{eq:lane_total_vcg_moment} \sum_{l \in \mathcal{L}} \sum_{c \in \mathcal{C}} \text{VCG}_l \cdot W_c \cdot H^C_c/2 \cdot x_{l,c} = m^V_{\text{cargo}} \\
    & \label{eq:lane_frame_start_lower_frame} 
    r_{c,f} \geq B_f \cdot \phi_{c,f} \quad & \forall c \in \mathcal{C}, f \in \mathcal{F} \\
    & \label{eq:lane_frame_start_lower_cargo} 
    r_{c,f} \geq (\text{lcg}_c - L^{\text{cargo}}_c/2) - M \cdot (1 - \phi_{c,f}) \quad & \forall c \in \mathcal{C}, f \in \mathcal{F} \\
    & \label{eq:lane_frame_end_upper_frame} 
    r_{c,f} + L^{\text{cargo}}_c \cdot \rho_{c,f} \leq B_{f+1} + M \cdot (1 - \phi_{c,f}) \quad & \forall c \in \mathcal{C}, f \in \mathcal{F}\setminus\{F\} \\
    & \label{eq:lane_frame_end_upper_cargo} 
    r_{c,f} + L^{\text{cargo}}_c \cdot \rho_{c,f} \leq (\text{lcg}_c + L^{C}_c/2) + M \cdot (1 - \phi_{c,f}) \quad & \forall c \in \mathcal{C}, f \in \mathcal{F} \\
    & \label{eq:lane_frame_start_upper} 
    r_{c,f} \leq M \cdot \phi_{c,f} \quad & \forall c \in \mathcal{C}, f \in \mathcal{F} \\
    & \label{eq:lane_frame_fraction_active} 
    \rho_{c,f} \leq \phi_{c,f} \quad & \forall c \in \mathcal{C}, f \in \mathcal{F} \\
    & \label{eq:lane_cargo_fraction_sum} 
    \sum_{f \in \mathcal{F}} \rho_{c,f} = \sum_{l \in \mathcal{L}} x_{l,c} \quad & \forall c \in \mathcal{C} \\
    & \label{eq:lane_frame_weight_distribution} 
    w_f = \sum_{c \in \mathcal{C}} W_c \cdot \rho_{c,f} \quad & \forall f \in \mathcal{F} \\
    &\label{eq:lane_domain_binary}
    x_{l,c},z_{l,c,c'},\phi_{c,f} \in \{0,1\} \quad & \forall l \in \mathcal{L}, c,c' \in \mathcal{C}, f\in \mathcal{F}\\   
    & \label{eq:lane_domain_merge} 
    \mu_p \in \{0,1\} \quad & \forall p \in \mathcal{M} \\
    & \label{eq:lane_domain_fraction} 
    \rho_{c,f} \in \mathbb{R} \quad & \forall c \in \mathcal{C}, f \in \mathcal{F} \\
    & \label{eq:lane_domain_position} 
    s_{l,c}, r_{c,f} \geq 0 \quad & \forall l \in \mathcal{L}, c \in \mathcal{C}, f \in \mathcal{F} \\
    & \label{eq:lane_domain_moments} 
    \text{lcg}_c, M_{lcg}, M_{tcg}, M_{vcg} \in \mathbb{R} \quad & \forall c \in \mathcal{C} \\
    & \text{Subject to stability constraints } \eqref{eq:stability_max_ballast_water} - \eqref{eq:stability_variable_domain_w} \nonumber
\end{align}

The objective function \eqref{eq:lane_objective} is similar to that of the hybrid and grid models. The lane formulation introduces two new variables for cargo placement. The binary variable $x_{l,c}$ controls if cargo $c$ is assigned to a lane $l\in \mathcal{L}$ where the set $\mathcal{L}$ contains all lanes, both single and merged lanes across the decks of the vessel. The variable $z_{l,c,c'}$ decides the precedence of cargo within a lane. Specifically, $z_{l,c,c'}$ determines if cargo $c$ is before cargo $c'$ in lane $l$. Constraints \eqref{eq:lane_precedence_required} establish the precedence relationship between cargo $c$ and $c'$ in lane $l$ by ensuring that if both cargo $c$ and $c'$ are placed in a lane then at least one of  $z_{l,c,c'}$ and  $z_{l,c',c}$ must be one. The constraints \eqref{eq:lane_precedence_exclusive} specify that the precedence between cargo $c$ and $c'$ can only be in a single direction. Constraints \eqref{eq:lane_precedence_bound_out} and \eqref{eq:lane_precedence_bound_in} ensure that if cargo $c$ is placed in lane $l$, then cargo $c$ will also have precedence in lane $l$. Further, it is an upper bound on the maximum precedence that can occur in lane $l$ through the big-M $M^{\max}_l$, which is the maximum number of cargo that can be assigned to a lane. From the precedence, the starting position of each cargo $s_{l,c'}$ can be found in constraints \eqref{eq:lane_no_overlap_forward}. The constraints ensure that if cargo $c'$ is placed right after cargo $c$, then the starting position $s_{l,c'}$ of $c'$ is after the starting position $s_{l,c}$ of $c$  plus the length of the cargo $L^{\text{cargo}}_c$. Note that the last part of the constraints $L^{\text{lane}}_l \cdot (1-z_{l,c,'})$ is a big-M to ensure that the starting position is only positive for $c,c'$ assigned to lane $l$. The restriction \eqref{eq:lane_position_bound} is an upper bound on the starting position $s_{l,c}$ that specifies that it cannot be greater than the length of the lane minus the length of the cargo. The width of the the lane and cargo is harder to account for than in the hybrid-slot and grid model, due to the the ,vessel's complex layout, and the TCG of a lane depends on its location on the deck. Constraints \eqref{eq:lane_merge_enable_new} ensure that if a merged lane is used, then $\mu_p$ must be active. The set $\mathcal{M}$ contains the merged lane $p$ and its parent lanes $(i,j)$. In constraints \eqref{eq:lane_merge_disable_original}, the original lanes become unavailable if $\mu_p$ is active. Constraints \eqref{eq:lane_width_constraint} ensure that the width of the cargo $c$ cannot exceed the width of the lane $l$. The LCG of cargo $c$ is found in constraints \eqref{eq:lane_cargo_lcg} from the starting position $s_{l,c}$, the length of the cargo $L^{\text{cargo}}_c$ and the start LCG of the lane $\text{LCG}^{\text{start}}_l$. From the individual cargo LCG, the total cargo LCG contribution is found in constraint \eqref{eq:lane_total_lcg_moment} where the LCG of the cargo is multiplied by its weight. The TCG contribution is found in a similar way in constraint \eqref{eq:lane_total_tcg_moment}, but the TCG is constant in the middle of the lane. In constraint \eqref{eq:lane_total_vcg_moment}, the VCG is found in the same manner as the TCG, but the height of the cargo is accounted for. 

The constraints \eqref{eq:lane_frame_start_lower_frame} to \eqref{eq:lane_frame_weight_distribution} determine which frames each cargo is present in. Again, this is more complex than in the other models since it cannot be pre-computed. Constraints \eqref{eq:lane_frame_start_lower_frame} introduce the variable $r_{c,f}$ as the starting position of cargo $c$ in frame $f$. $r_{c,f}$ must be greater than either $B_f$ which is the start position of frame $f$ or the starting position of cargo $c$ in \eqref{eq:lane_frame_start_lower_cargo}. The binary variable $\phi_{c,f}$ is introduced to ensure that $r_{c,f}$ is zero when cargo $c$ is outside a frame. Constraints \eqref{eq:lane_frame_end_upper_frame} and \eqref{eq:lane_frame_end_upper_cargo} do the same for bounding the ending position. Here, the variable $p_{f,c}$ is introduced as the fraction of how much cargo $c$ is in frame $f$. Constraints \eqref{eq:lane_frame_start_upper} and \eqref{eq:lane_frame_fraction_active} ensure that $r_{c,f}$ and $p_{f,c}$ are active when cargo is present in a frame. Constraints \eqref{eq:lane_cargo_fraction_sum} specify that if cargo is placed, then its fractions in all frames must sum to 1. Finally, in constraints \eqref{eq:lane_frame_weight_distribution}, the total cargo weight in each frame is determined.


\section{Computational Study} \label{sec:results}
In collaboration with our industrial partner, a large Danish RoRo shipping company, we have designed a set of benchmark instances to study the three model representations. We consider three vessel classes: small, medium, and large. Vessel stability data, such as hydrostatic data tables, have been provided by our industrial partner from actual vessel. Though the instances are confidential, we strive to share details. The number of decks, the nominal trailer capacity, maximum tonnage, and number of ballast tanks are reported in the \autoref{tab:instance_data_vessel}.

\begin{table}[!thb]
\footnotesize
\centering
\begin{tabular}{lrrrc}
\hline
\multicolumn{5}{|c|}{\textbf{VESSEL SPECIFICATIONS}}  \\ \hline
\multicolumn{1}{|l|}{\textbf{Size}}           & \multicolumn{1}{c|}{\textbf{Light Ship Weight(T)}} & \multicolumn{1}{c|}{\textbf{Decks}}    & \multicolumn{1}{c|}{\textbf{Trailers}}        & \multicolumn{1}{c|}{\textbf{Ballast Tanks}} \\ \hline
\multicolumn{1}{|l|}{Small}                   & \multicolumn{1}{r|}{5773}                          & \multicolumn{1}{r|}{3}                 & \multicolumn{1}{r|}{113}                      & \multicolumn{1}{r|}{21}                     \\ \hline
\multicolumn{1}{|l|}{Medium}                  & \multicolumn{1}{r|}{10443}                         & \multicolumn{1}{r|}{4}                 & \multicolumn{1}{r|}{262}                      & \multicolumn{1}{r|}{22}                     \\ \hline
\multicolumn{1}{|l|}{Large}                   & \multicolumn{1}{r|}{18126}                         & \multicolumn{1}{r|}{5}                 & \multicolumn{1}{r|}{428}                      & \multicolumn{1}{r|}{25}                     \\ \hline
\multicolumn{1}{c}{}                          & \multicolumn{1}{c}{}                               & \multicolumn{1}{c}{}                   & \multicolumn{1}{c}{}                          & 
\end{tabular}
\caption{Specification of vessel data.}
\label{tab:instance_data_vessel}
\end{table}

The instances include a set of heterogeneous cargo composed of four main types: trailers, cars, machinery (e.g., agricultural vehicles), and SECU boxes (wheeled platforms for the transport of breakbulk or containers). The size and weight range of each type can be found in \autoref{tab:instance_data_cargo}.

\begin{table}[!thb]
\footnotesize
\centering
\begin{tabular}{lrrr}
\hline
\multicolumn{4}{|c|}{\textbf{CARGO TYPE SPECIFICATIONS}}                                                                                                               \\ \hline
\multicolumn{1}{|l|}{\textbf{Type}} & \multicolumn{1}{c|}{\textbf{Length(m)}} & \multicolumn{1}{c|}{\textbf{Width(m)}} & \multicolumn{1}{c|}{\textbf{Weight Range(T)}} \\ \hline
\multicolumn{1}{|l|}{Trailer}       & \multicolumn{1}{r|}{13.6}               & \multicolumn{1}{r|}{2.5}               & \multicolumn{1}{r|}{7-30}                     \\ \hline
\multicolumn{1}{|l|}{Car}           & \multicolumn{1}{r|}{4.0}                & \multicolumn{1}{r|}{1.5}               & \multicolumn{1}{r|}{1-3}                      \\ \hline
\multicolumn{1}{|l|}{Machinery}     & \multicolumn{1}{r|}{4.6}                & \multicolumn{1}{r|}{4.5}               & \multicolumn{1}{r|}{5-10}                     \\ \hline
\multicolumn{1}{|l|}{SECU}          & \multicolumn{1}{r|}{13.8}               & \multicolumn{1}{r|}{3.6}               & \multicolumn{1}{r|}{40-80}                    \\ \hline
\multicolumn{1}{c}{}                & \multicolumn{1}{c}{}                    & \multicolumn{1}{c}{}                   & \multicolumn{1}{c}{}                         
\end{tabular}
\caption{Specification of cargo types.}
\label{tab:instance_data_cargo}
\end{table}

To test the models across different shipping services, we have devised 6 cargo-type distributions. \autoref{tab:instance_data_cargo_dist} shows the percentage of cargo from the booking list that falls into each of the cargo categories. The distributions have been designed both to resemble what our industrial partner sees in practice and to test the capability of the models. A clear example of the latter is a distribution in which only trailers are used, or in which cars are not transported.

\begin{table}[!thb]
\footnotesize
\centering
\begin{tabular}{lrrrc}
\hline
\multicolumn{5}{|c|}{\textbf{CARGO TYPE DISTRIBUTIONS (\%)}}                                                                                                                                                                              \\ \hline
\multicolumn{1}{|l|}{\textbf{Scenario}}       & \multicolumn{1}{c|}{\textbf{Trailers}}             & \multicolumn{1}{c|}{\textbf{Cars}}     & \multicolumn{1}{c|}{\textbf{Machinery}}       & \multicolumn{1}{c|}{\textbf{SECU}}     \\ \hline
\multicolumn{1}{|l|}{Equal Distribution (ED)} & \multicolumn{1}{r|}{25}                            & \multicolumn{1}{r|}{25}                & \multicolumn{1}{r|}{25}                       & \multicolumn{1}{r|}{25}                     \\ \hline
\multicolumn{1}{|l|}{Mixed Load (ML)}         & \multicolumn{1}{r|}{50}                            & \multicolumn{1}{r|}{30}                & \multicolumn{1}{r|}{10}                       & \multicolumn{1}{r|}{10}                     \\ \hline
\multicolumn{1}{|l|}{Trailer Only (TO)}       & \multicolumn{1}{r|}{100}                           & \multicolumn{1}{r|}{0}                 & \multicolumn{1}{r|}{0}                        & \multicolumn{1}{r|}{0}                      \\ \hline
\multicolumn{1}{|l|}{Trailer \& SECU (TS)}    & \multicolumn{1}{r|}{75}                            & \multicolumn{1}{r|}{0}                 & \multicolumn{1}{r|}{0}                        & \multicolumn{1}{r|}{25}                     \\ \hline
\multicolumn{1}{|l|}{Trailer Heavy (TH)}      & \multicolumn{1}{r|}{80}                            & \multicolumn{1}{r|}{6.67}              & \multicolumn{1}{r|}{6.67}                     & \multicolumn{1}{r|}{6.67}                   \\ \hline
\multicolumn{1}{|l|}{No Cars (NC)}            & \multicolumn{1}{r|}{70}                            & \multicolumn{1}{r|}{0}                 & \multicolumn{1}{r|}{15}                       & \multicolumn{1}{r|}{15}                     \\ \hline
\end{tabular}
\caption{Specification of cargo distributions.}
\label{tab:instance_data_cargo_dist}
\end{table}

Each deck area is converted to match the models' representation. For the hybrid-slot models, the industrial partner has provided us with a set of pre-computed slot positions for standard 14x2.5-meter trailers. Slots of other cargo sizes are computed based on those. For the lane model, the maximum number of lanes is computed assuming a minimum lane width of 2.5 meters. Finally, for the grid model, each grid cell is 2.5x2.5 meters. This grid size is too coarse; a finer mesh would certainly yield better results. However, it was not possible to solve any instances (within reasonable computational times) with the grid model for any grid finer than the one used in this study. 

A total of 162 instances have been generated, and \autoref{tab:instance_component} summarises their components. For each combination of cargo distribution and vessel size, 6 instances are generated, distinguished by the booked capacity and the sampling method for the cargo weight. Booking capacity refers to the percentage of the nominal vessel capacity that is present in the cargo booking list. We operate with 3 different cases: under-booking (80\%), full-booking (100\%), and over-booking (120\%). The weights of the cargo are sampled from their relative ranges ($[w_{min}, w_{max}]$) via a truncated normal distribution  $\mathcal{N}(\mu, \sigma)$. The cargo weights are generated using a truncated normal distribution $\mathcal{N}(\mu, \sigma)$ bounded by $[w_{min}, w_{max}]$, where $\mu$ is set at 20\%, 50\%, or 80\% of the weight range for light, medium, and heavy profiles respectively, and $\sigma = (w_{max} - w_{min})/6$. Each instance is named according to the format \textit{vessel\_type\_ratio\_weight}, for example, \textit{S\_ED\_60\_L} or \textit{L\_NC\_120\_H}, where \textit{vessel} denotes the vessel size (S for small, M for medium, L for large), \textit{type} represents the cargo distribution type (e.g., ED for Equal Distribution, NC for No-Cars), \textit{ratio} indicates the percentage of the vessel's capacity used (e.g. 80, 100, 120), and \textit{weight} specifies the cargo weight category (L for light, M for medium, H for heavy). 

\begin{table}[!t]
\centering
\begin{tabular}{|lrl|}
\hline
\multicolumn{3}{|c|}{\textbf{Instance Overview}}                                                                         \\ \hline
\multicolumn{1}{|c|}{\textbf{Component}}       & \multicolumn{1}{c|}{\textbf{Count}} & \multicolumn{1}{c|}{\textbf{Description}} \\ \hline
\multicolumn{1}{|l|}{Cargo Distribution}             & \multicolumn{1}{r|}{6}              & ED, ML, TO, TS, TH, NC                    \\ \hline
\multicolumn{1}{|l|}{Capacity booked (\%)}              & \multicolumn{1}{r|}{3}              & 80, 100, 120                        \\ \hline
\multicolumn{1}{|l|}{Weight Categories}        & \multicolumn{1}{r|}{3}              & Light, Medium, Heavy                      \\ \hline
\multicolumn{1}{|l|}{Ship Size}                & \multicolumn{1}{r|}{3}              & Small, Medium, Large                      \\ \hline
\multicolumn{1}{|l|}{\textbf{Total Instances}} & \multicolumn{1}{r|}{\textbf{162}}   & \textbf{6 × 3 × 3 × 3}                    \\ \hline
\end{tabular}
\caption{The four components defining an instance.}\label{tab:instance_component}
\end{table}

\subsection{Results}
The computational study is implemented in Julia using the JuMP modelling package, and all problem instances are solved with Gurobi 12.0.0. All experiments are conducted on an Intel\textsuperscript{\textregistered} Xeon\textsuperscript{\textregistered} Gold 6226R processor, with 8 available threads. We impose a 1-hour time limit per instance. The models are evaluated on the following criteria: the number of instances with a solution, the number of instances solved to proven optimality, the average optimality gap percentage reported by Gurobi, the percentage of cargo successfully allocated from the booking list, and the time to solve them. The aggregated results are presented in \autoref{tab:grand_summary_final}, calculated as the average performance across all instances for which a feasible solution was obtained. Detailed, instance-specific outcomes are provided in the appendix. The models are evaluated under three configurations: 
\begin{description}
    \item[No Stab:] No stability constraints are imposed; the objective is to maximise the total cargo loaded, thereby assessing the model’s raw placement capability.
    \item[Simple Stab:]  Stability constraints related to centres of gravity and ballast tank adjustments are applied.
    \item [Full Stab:] Full seaworthiness conditions, including stress force constraints.
\end{description}
Together, these three scenarios assess the models’ performance and illustrate the impact of stability considerations on their overall effectiveness.

\begin{table}[!htb]
\centering
\small 
\setlength{\tabcolsep}{2pt} 
\renewcommand{\arraystretch}{1.2} 

\caption{Aggregated performance comparison including Cargo Utilization. 
\textnormal{Values represent averages. \textbf{Util}: Cargo Utilization (\%), \textbf{Gap}: Optimality Gap (\%), \textbf{Time}: Time (s), \textbf{O/S}: Optimal/Solved instances. There are 54 instances per vessel size and 162 in total.}}
\label{tab:grand_summary_final}

\begin{tabular}{@{}ll cccc cccc cccc@{}}
\toprule
& & \multicolumn{4}{c}{\textbf{Hybrid}} & \multicolumn{4}{c}{\textbf{Grid}} & \multicolumn{4}{c}{\textbf{Lane}} \\
\cmidrule(lr){3-6} \cmidrule(lr){7-10} \cmidrule(lr){11-14}
\textbf{Mode} & \textbf{Size} & \textbf{Util} & \textbf{Gap} & \textbf{Time} & \textbf{O/S} & \textbf{Util} & \textbf{Gap} & \textbf{Time} & \textbf{O/S} & \textbf{Util} & \textbf{Gap} & \textbf{Time} & \textbf{O/S} \\
\midrule

\multirow{4}{*}{\shortstack[l]{\textbf{No}\\ \textbf{Stab.}}} 
 & Small  & 97.11 & 0.02 & 325 & 50/54 & 89.83 & 0.00 & 10 & 54/54 & 69.28 & 180.83 & 2836 & 12/51 \\
 & Medium & 95.65 & 0.33 & 833 & 45/54 & 90.44 & 0.00 & 69 & 54/54 & 64.36 & 107.17 & 3600 & 0/54 \\
 & Large  & 96.58 & 0.58 & 855 & 39/48 & 91.02 & 0.01 & 211 & 50/54 & --- & --- & --- & 0/0 \\
 \cmidrule(l){2-14} 
 & \textit{Avg} & \textbf{96.45} & \textit{0.31} & \textit{671} & \textit{134/156} & \textit{89.10} & \textbf{0.01} & \textbf{211} & \textbf{158/162} & \textit{67.31} & \textit{151.37} & \textit{3143} & \textit{12/79} \\
\midrule

\multirow{4}{*}{\shortstack[l]{\textbf{Simple}\\ \textbf{Stab.}}} 
 & Small  & 97.10 & 0.16 & 1001 & 41/54 & 85.21 & 0.02 & 529 & 50/54 & 68.73 & 121.11 & 2906 & 11/43 \\
 & Medium & 95.08 & 0.89 & 1421 & 33/53 & 89.18 & 0.07 & 1939 & 25/53 & 61.17 & 78.06 & 3454 & 0/8 \\
 & Large  & 94.61 & 3.46 & 1408 & 36/54 & 91.76 & 0.03 & 1734 & 30/49 & --- & --- & --- & 0/0 \\
 \cmidrule(l){2-14} 
 & \textit{Avg} & \textbf{95.60} & \textit{1.51} & \textbf{1277} & \textbf{110/161} & \textit{88.72} & \textbf{0.04} & \textit{1401} & \textit{105/156} & \textit{66.21} & \textit{106.76} & \textit{3125} & \textit{11/51} \\
\midrule

\multirow{4}{*}{\shortstack[l]{\textbf{Full}\\ \textbf{Stab.}}} 
 & Small  & 97.10 & 0.08 & 945 & 43/54 & 84.53 & 0.02 & 569 & 50/54 & 91.58 & 16.87 & 3405 & 7/12 \\
 & Medium & 95.08 & 0.79 & 1381 & 34/53 & 89.03 & 0.07 & 1950 & 25/52 & --- & --- & --- & 0/0 \\
 & Large  & 93.06 & 6.18 & 1515 & 33/43 & 89.97 & 0.03 & 1883 & 29/53 & --- & --- & --- & 0/0 \\
 \cmidrule(l){2-14} 
 & \textit{Avg} & \textbf{95.08} & \textit{2.35} & \textbf{1280} & \textbf{110/160} & \textit{87.84} & \textbf{0.04} & \textit{1467} & \textit{104/159} & \textit{---} & \textit{---} & \textit{---} & \textit{7/12} \\

\bottomrule
\end{tabular}
\end{table}

Overall, the hybrid-slot and grid models perform well, whereas the lane model performs poorly. The lane model can find solutions for the small and medium-sized vessels, but fails to find feasible solutions for the large vessel, exhibiting an average gap between 60-70\% and reaching the time limit of one hour. The comparison between the hybrid-slot and grid models is more nuanced. For the small vessel, the grid model outperforms in all metrics except cargo utilisation. For the medium and large vessels, this trend reverses, and the grid model is only better in terms of the gap. The large discrepancy in cargo utilisation between the hybrid-slot and grid model is due to the grid resolution of 2.5x2.5m. This leads to suboptimal space efficiency, especially when placing cars, as they are only 1.5m wide, resulting in wasted space. \autoref{tab:grid_low_res} shows that the grid model without stability outperforms the hybrid-slot when using a grid resolution of 0.5x0.5m at the cost of increasing the computational time by almost a factor of 4. For the hybrid-slot model, adding stability constraints doubles the average computational time and decreases the cargo loaded by around 1\%. The same pattern is observed with the grid model, but the average solution time increases by almost 6.5 times, with no significant difference when stress forces are added. The lane model is difficult to judge, as so few instances are solved to optimality, but performance degrades significantly when stability is added. Based on our experiments, the hybrid-slot model achieves the best overall performance, followed closely by the grid model. The results align with prior research, which indicates that grid-based and lane-based formulations require specialised heuristics to yield good solutions on realistic data \cite{ovstebo_optimization_2011,hansen_adaptive_2020}.

A granular analysis of instance-specific performance reveals distinct behaviours across the three formulations. While the hybrid-slot model is generally robust, its efficiency degrades when the cargo diversity increases, particularly in the Equally Distributed (ED) and Mixed Load (ML) instances (\autoref{tab:detailed_breakdown_hybrid}). In these scenarios, there is a sharp increase in computational time and a corresponding decrease in the number of proven optimal solutions. The grid model exhibits a similar sensitivity to the cargo mix but is notably more affected by vessel size (\autoref{tab:detailed_breakdown_grid}). This is expected, as the combinatorial complexity of the grid formulation scales more aggressively with vessel dimensions than the hybrid-slot model. Finally, the lane models' performance degrades drastically when instances contain cars (\autoref{tab:detailed_breakdown_grid}). This degradation is attributable to the precedence constraints, which grow rapidly in number when smaller units, such as cars, are introduced to the booking list.

The results show that Gurobi cannot solve all instances to optimality. From an industry perspective, optimality is not strictly necessary, but the time to get good solutions is critical. Since the average solution time for both the hybrid-slot and grid models exceeds 15 minutes, it could be argued that these are unsuitable for actual use cases, as planning in a terminal is a dynamic, reactive environment. Therefore, there remains a significant opportunity to develop good heuristic for the RSSP. From a practical perspective, the lane model does not appear to be viable for stowage planning involving stability calculations. The comparison between the grid and hybrid slot is more nuanced. The grid model does produce better solutions when the grid size is small enough. Further, the grid model performs better when there are many different cargo types present. However, in reality, when transporting trailers on RoRo vessels, so called "SAT trestles" are used to secure the trailers. SAT trestles are metalling support structures that, attached to a trailer, can lock into the deck holes. These holes are spaced according to a predefined pattern, making them ideal for slot representation. There seems not to be a single superior approach. Depending on the problem setting, the grid or hybrid-slot model might be preferable.

\begin{table}[!htb]
\centering
\small
\begin{threeparttable}
    \caption{Comparative Performance of Hybrid vs. Grid Models with Ballast. 
    \textnormal{\footnotesize Comparison across Simple Stability and All Stability modes. Time in seconds, Ballast in tonnes.}}
    \label{tab:hybrid_vs_grid_ballast}

    \begin{tabular}{@{}ll rrr rrr@{}}
    \toprule
    & & \multicolumn{3}{c}{\textbf{Simple Stability}} & \multicolumn{3}{c}{\textbf{Full Stability}} \\
    \cmidrule(lr){3-5} \cmidrule(lr){6-8}
    \textbf{Ship} & \textbf{Model} & \textbf{Cargo} & \textbf{Ballast} & \textbf{Time} & \textbf{Cargo} & \textbf{Ballast} & \textbf{Time} \\
     & & \textbf{(\%)} & \textbf{(t)} & \textbf{(s)} & \textbf{(\%)} & \textbf{(t)} & \textbf{(s)} \\
    \midrule

    \multirow{2}{*}{\textbf{Small}} 
    & Hybrid & 97.10 & 1745 & 1001.24 & 97.10 & 1861 & 944.55 \\
    & Grid   & 85.21          & 2341          & 529.47  & 84.53          & 2445          & 569.21 \\
    \addlinespace

    \multirow{2}{*}{\textbf{Medium}} 
    & Hybrid & 95.08 & 2326 & 1420.80 & 95.08 & 2559 & 1380.97 \\
    & Grid   & 89.18          & 3000          & 1938.63          & 89.03          & 3071          & 1950.20 \\
    \addlinespace

    \multirow{2}{*}{\textbf{Large}} 
    & Hybrid & 94.61 & 3468 & 1407.56 & 93.06 & 3951 & 1514.67 \\
    & Grid   & 91.76          & 4592          & 1733.70          & 89.97          & 5325          & 1882.63 \\

    \midrule
    \multicolumn{2}{@{}l}{\textbf{Grand Total}} & 92.16 & 2912 & 1338.57 & 91.46 & 3202 & 1373.71 \\
    \bottomrule
    \end{tabular}
\end{threeparttable}
\end{table}

\begin{figure}[!htb]
     \centering
     \begin{subfigure}[t]{0.45\textwidth}
         \centering
         \includegraphics[width=\textwidth]{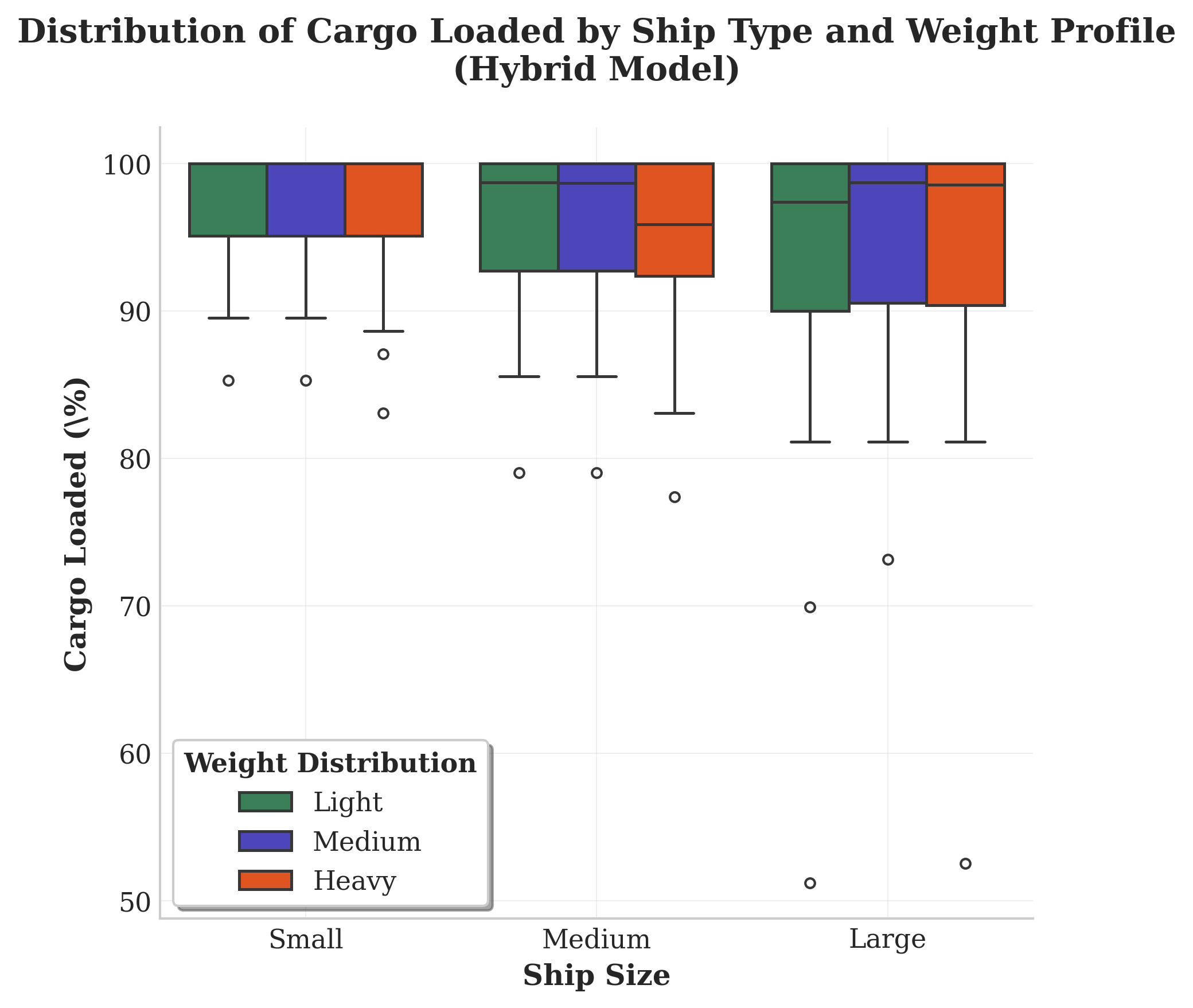}
         \caption{The hybrid models booking utilization for each weight distribution and ship size. The plot highlights that the weight distribution does not significantly impact the booking utilization.}
         \label{fig:booking_weight}
     \end{subfigure}
     \hfill
     \begin{subfigure}[t]{0.45\textwidth}
         \centering
         \includegraphics[width=\textwidth]{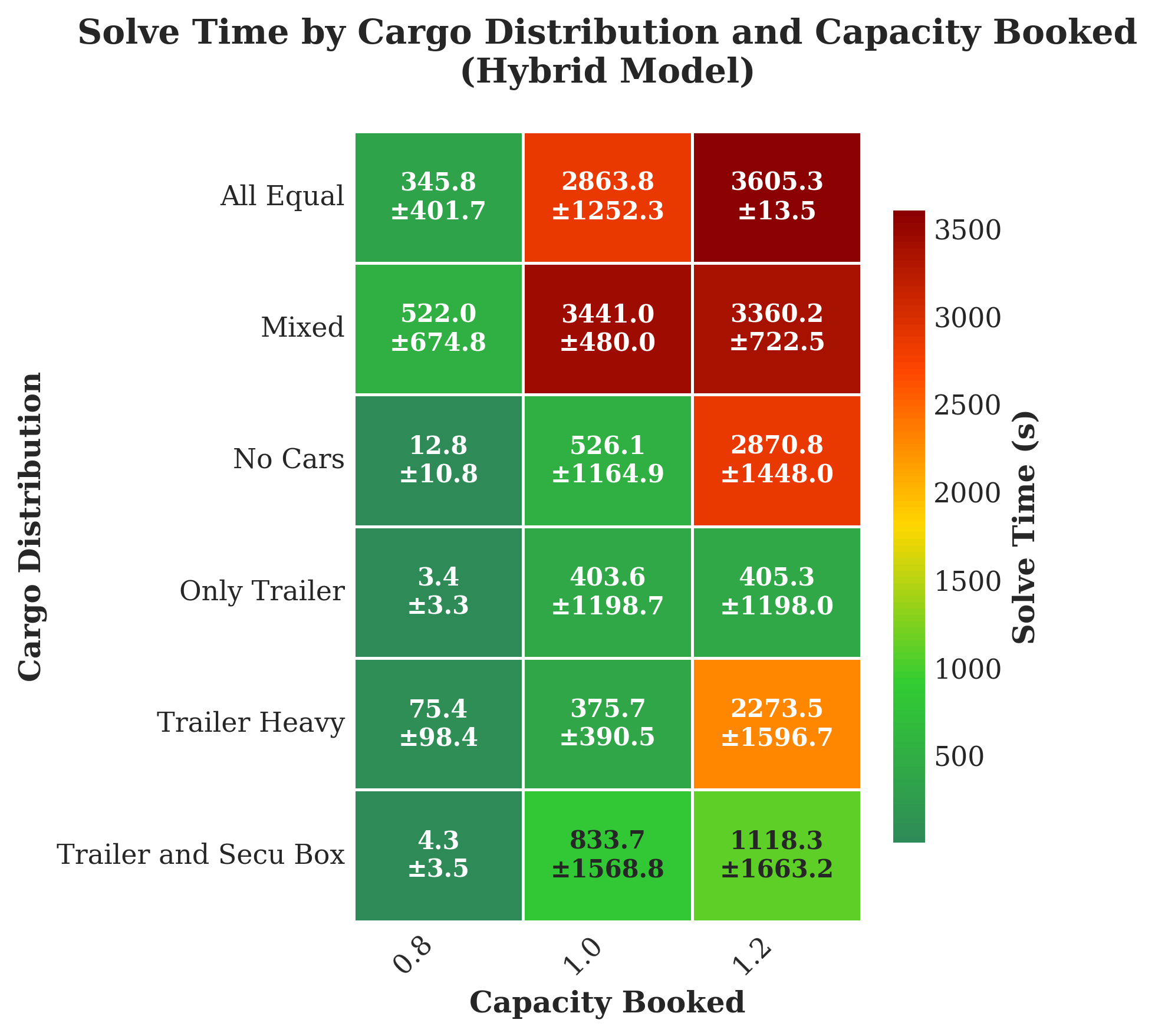}
         \caption{Overview of the cargo distribution and capacity solution time for the hybrid model. The numbers show the average solution time (s) and standard deviation. The plot highlights the relationship between the capacity booked, cargo type distribution and the solution time.}
         \label{fig:ballast_weight}
     \end{subfigure}
    \label{fig:boxplots}
    \caption{Overview of the metrics of the hybrid-slot model}    
\end{figure}

\autoref{tab:hybrid_vs_grid_ballast} shows a negligible difference in solution times between the Simple and Full stability configurations, indicating that the stress constraints are not a computational bottleneck. However, there is a significant increase in average ballast water utilisation when moving from Simple to Full stability. This disparity demonstrates that the resulting stowage plans differ and should therefore be added to ensure true seaworthiness. \autoref{fig:booking_weight} shows the cargo loaded for each vessel and the weight distribution for the hybrid-slot model. The plot shows that there are no significant differences in the amount of cargo loaded across the different weight distributions, indicating no positive correlation. \autoref{fig:ballast_weight} shows the execution time of the hybrid-slot model's capacity booked and cargo distribution. The plot shows that all vessels solve the cargo type distributions \textit{TO} and \textit{TC} quickly. The other cargo type distributions are slower to solve, and their performance degrades significantly at high-capacity bookings.

\subsection{Hazardous and Refrigerated Extension} \label{sec:extension}
A second aspect of RoRo stowage planning is how to accommodate refrigerated and hazardous cargo. Refrigerated cargo must be placed near electrical outlets on deck. The planning for hazardous cargo involves strategically placing these items in accordance with the stowage restrictions outlined in the International Maritime Dangerous Goods (IMDG) Code \cite{international_maritime_organization_international_2014}. For instance, a cargo classified as IMDG class 4.1 (flammable solid) must be segregated by a distance of 6 meters from a cargo classified as class 6.1 (toxic substance), adhering to a segregation degree of 2. In this study, we consider 18 distinct classes of hazardous materials and five segregation rules, ranging from 3 meters to 36 meters of separation. We only consider segregation between cargo and non-crew areas, engine areas, or specific required locations like "cargo must be on the weather deck". However, this could be implicitly modelled by restricting the set of slots or positions where cargo can be placed. As one might anticipate, planning hazardous cargo becomes more complex when a significant amount is involved in a voyage. To analyse the impact of specialised cargo on the optimisation performance, we have extended the hybrid-slot model with the following constraints to accommodate refrigerated and hazardous cargo.

\begin{align}
    & \label{eq:hybrid_ref_slot} x_{c,s} = 0 & \quad \forall c \in C^R, s \in S \setminus S^R \\
    & \label{eq:hybrid_haz_restriction} x_{c,s} + x_{c',s'} \leq 1 & \quad  \forall  (c,s),(c',s') \in \mathcal{R}
\end{align}

Constraints \eqref{eq:hybrid_ref_slot} ensure that refrigerated cargo is allocated to a slot position that can reach a power outlet. Constraints \eqref{eq:hybrid_haz_restriction} ensure minimal segregation distance between conflicting hazardous cargo. The set $\mathcal{R}$ contains the combinations of cargo locations $x_{c,s}$ that cannot be utilised concurrently by $c$ and $c'$ due to their hazardous content. The rules are obtained from the International Maritime Organisation \cite{international_maritime_organization_international_2014}. 

The quantity of hazardous segregation constraints is considerable, as they are present between all slots and hazardous cargo. This becomes particularly challenging when the number of hazardous cargo units increases, resulting in greater complexity and making the model harder to build and solve. To address this issue, we introduce lazy constraints. Unlike traditional constraints, lazy constraints are not added during the initial model construction. Instead, they are introduced only during the search procedure if a solution violates these constraints. Lazy constraints are particularly effective when they are not always necessary. Given that the proportion of hazardous cargo on a voyage is typically low, constituting only 5-10\% of the total cargo, lazy constraints are an ideal approach, as it is expected that not all of them will need to be added to the solution. If a violation is detected, a new constraint is introduced for the specific cargo slots, formulated as  $x_{c,s} + x_{c',s'} \leq 1$. 

We test the refrigerated and hazardous extensions on instances with no cars, medium weight, and a 100\% booking list length for all vessels. Refrigeration ability is randomly assigned to half of the available slots. Hazardous cargo classes are incorporated into the instance, with two distinct generation methods employed: Synthetic and Realistic. The Synthetic generation method samples uniformly from the 17 different hazardous cargo classes, representing an extreme rather than a realistic scenario. In contrast, the Realistic generation method samples from a distribution that reflects real-world conditions, where only a few classes have high probabilities of occurrence. Finally, we generate instances with varying proportions of refrigerated and hazardous cargo, ranging from 5\% to 25\%.

\begin{table}[H]
\centering
\small
\begin{threeparttable}
    \caption{Overview of Refrigerated and Hazardous Results. 
    \textnormal{\footnotesize The results are an average across all three vessels. Refrigerated-only results are excluded as they are trivially solved.}}
    \label{tab:haz_ref_results}

    \begin{tabular}{@{}lll rr cc@{}}
    \toprule
    \multirow{2}{*}{\textbf{Cargo Type}} & \multirow{2}{*}{\textbf{Ratio}} & \multirow{2}{*}{\textbf{Distribution}} & \textbf{Time} & \textbf{Gap} & \multicolumn{2}{c}{\textbf{Utilization (\%)}} \\
    \cmidrule(lr){6-7}
     & & & \textbf{(s)} & \textbf{(\%)} & \textbf{Haz} & \textbf{Ref} \\
    \midrule

    \multirow{8}{*}{\textbf{Hazardous}} 
     & \multirow{2}{*}{5\%}  & Synthetic & 1347.70 & 83.23 & 74 & 0 \\
     &                       & Realistic & 124.79  & 0.00  & 100 & 0 \\
     \addlinespace 
     & \multirow{2}{*}{10\%} & Synthetic & 3621.89 & 44.81 & 59 & 0 \\
     &                       & Realistic & 158.11  & 0.00  & 100 & 0 \\
     \addlinespace
     & \multirow{2}{*}{15\%} & Synthetic & 1801.58 & 46.91 & 76 & 0 \\
     &                       & Realistic & 744.12  & 0.00  & 100 & 0 \\
     \addlinespace
     & \multirow{2}{*}{25\%} & Synthetic & 3604.52 & 132.49 & 48 & 0 \\
     &                       & Realistic & 2302.56 & 0.01   & 100 & 0 \\
    \midrule

    \multirow{8}{*}{\shortstack[l]{\textbf{Hazardous \&}\\\textbf{Refrigerated}}} 
     & \multirow{2}{*}{5\%}  & Synthetic & 3610.93 & 2.48  & 82 & 100 \\
     &                       & Realistic & 275.66  & 0.00  & 100 & 100 \\
     \addlinespace
     & \multirow{2}{*}{10\%} & Synthetic & 3600.15 & 6.88  & 88 & 96 \\
     &                       & Realistic & 282.73  & 0.00  & 100 & 100 \\
     \addlinespace
     & \multirow{2}{*}{15\%} & Synthetic & 3604.84 & 85.20 & 67 & 79 \\
     &                       & Realistic & 733.10  & 0.00  & 100 & 100 \\
     \addlinespace
     & \multirow{2}{*}{25\%} & Synthetic & 3618.09 & 227.87 & 49 & 67 \\
     &                       & Realistic & 1174.53 & 0.00   & 100 & 100 \\

    \bottomrule
    \end{tabular}
\end{threeparttable}
\end{table}

\autoref{tab:haz_ref_results} shows that the realistic hazardous instances can easily be solved to optimality, and that all hazardous cargo is loaded. Instances involving only refrigerated cargo are not shown, as they are also easily resolved. Furthermore, the combined refrigerated and hazardous instances do not appear to be significantly impacted by the presence of refrigerated cargo. This indicates that the number of refrigerated slots is not the bottleneck. In contrast, the Synthetic instances are shown to be harder. Many instances exhibit high computational times and fail to place all hazardous cargo. These results highlight the sensitivity of the model towards an increasing number of hazardous cargo classes. 

\subsection{Stability remarks for the hybrid-slot model}
Modelling the complex stability equations requires several assumptions and simplifications. The hydrostatic data is discretised, the centres of gravity of the ballast tanks are linearised, and the total number of frames is reduced. To assess whether these modelling choices are acceptable, we validate all solutions from the hybrid-slot model with the vessels’ stability software, Loadmaster\textsuperscript{\textcopyright} X5 (Kockumation), which performs a full stability and strength check for each loading condition. As Loadmaster\textsuperscript{\textcopyright} X5 is proprietary software for checking vessel sea worthiness. The software can be used for manual planning, but has not automated features. Confidentiality hinder us from disclosing the exact evaluation metrics used by the software; instead, we show the number of stability errors our solutions would yield if used. \autoref{tab:detailed_stability_errors} summarises the stability measures that generate errors in Loadmaster for the tested hybrid-slot model instances.

\begin{table}[!htb]
\centering
\scriptsize 
\setlength{\tabcolsep}{3pt} 
\renewcommand{\arraystretch}{1.1}

\begin{threeparttable}
    \caption{Detailed Breakdown of Stability Errors by Cargo Configuration. 
    \textnormal{View of instances violating stability limits.}}
    \label{tab:detailed_stability_errors}

    \begin{tabular}{@{}ll ll r c r@{}}
    \toprule
    \textbf{Vessel} & \textbf{Config} & \textbf{Distribution} & \textbf{Measure} & \textbf{Err.} & \textbf{Mean} & \textbf{Limit} \\
    \midrule

    \multirow{17}{*}{\textbf{Small}} 
    & \multirow{5}{*}{\textbf{Heavy 100}} 
      & No Cars & GM Actual & 1 & 1.12\,m & --- \\
    & & \multirow{2}{*}{Only Trailer} & GM Actual & 1 & 1.27\,m & --- \\
    & & & Prop. Imm. & 1 & 47.0\% & $<50$\% \\
    & & Trailer + Secu & GM Actual & 1 & 1.13\,m & --- \\
    \cmidrule{2-7}

    & \multirow{7}{*}{\textbf{Heavy 120}} 
      & \multirow{2}{*}{All Equal} & GM Actual & 1 & 1.15\,m & --- \\
    & & & Prop. Imm. & 1 & 44.0\% & $<50$\% \\
    & & \multirow{2}{*}{Mixed} & GM Min & 1 & --- & --- \\
    & & & Prop. Imm. & 1 & 39.0\% & $<50$\% \\
    & & No Cars & GM Actual & 1 & 1.06\,m & --- \\
    & & Only Trailer & Prop. Imm. & 1 & 49.0\% & $<50$\% \\
    & & Trailer + Secu & GM Actual & 1 & 1.17\,m & --- \\
    \cmidrule{2-7}

    & \multirow{5}{*}{\textbf{Heavy 80}} 
      & \multirow{2}{*}{No Cars} & GM Min & 1 & --- & --- \\
    & & & Heel & 1 & 1.10° & $<1.0$° \\
    & & Only Trailer & Prop. Imm. & 1 & 44.0\% & $<50$\% \\
    & & Trailer + Secu & GM Actual & 1 & 1.20\,m & --- \\
    & & Trailer Heavy & Prop. Imm. & 1 & 44.0\% & $<50$\% \\
    \midrule

    \multirow{13}{*}{\textbf{Medium}} 
    & \multirow{3}{*}{\textbf{Heavy 100}} 
      & \multirow{2}{*}{All Equal} & GM Actual & 1 & 0.96\,m & --- \\
    & & & Heel & 1 & 2.30° & $<2.0$° \\
    & & No Cars & GM Actual & 1 & 0.96\,m & --- \\
    \cmidrule{2-7}

    & \multirow{3}{*}{\textbf{Heavy 80}} 
      & All Equal & Heel & 1 & 2.10° & $<2.0$° \\
    & & \multirow{2}{*}{Trailer + Secu} & GM Actual & 1 & 1.00\,m & --- \\
    & & & Heel & 1 & 2.10° & $<2.0$° \\
    \cmidrule{2-7}

    & \multirow{7}{*}{\textbf{Med. 80}} 
      & \multirow{2}{*}{Mixed} & GM Min & 1 & --- & --- \\
    & & & Trim & 1 & -0.97\,m & \tiny{$-0.95, 1.5$} \\
    & & No Cars & Trim & 1 & -0.97\,m & \tiny{$-0.95, 1.5$} \\
    & & \multirow{2}{*}{Trailer + Secu} & GM Min & 1 & --- & --- \\
    & & & Trim & 1 & -0.99\,m & \tiny{$-0.95, 1.5$} \\
    & & \multirow{2}{*}{Trailer Heavy} & GM Min & 1 & --- & --- \\
    & & & Trim & 1 & -0.97\,m & \tiny{$-0.95, 1.5$} \\
    \midrule

    \multirow{4}{*}{\textbf{Large}} 
    & \multirow{2}{*}{\textbf{Heavy 100}} 
      & No Cars & KG & 1 & 15.65\,m & $<15.57$ \\
    & & Trailer + Secu & KG & 1 & 15.83\,m & $<15.57$ \\
    \cmidrule{2-7}

    & \multirow{2}{*}{\textbf{Heavy 120}} 
      & No Cars & KG & 1 & 15.64\,m & $<15.57$ \\
    & & Trailer + Secu & KG & 1 & 15.64\,m & $<15.57$ \\

    \bottomrule
    \end{tabular}
    
    \begin{tablenotes}
        \scriptsize
        \item \textbf{Notes:} Prop. Imm. = Propeller Immersion.
        \item Med. 80 = Medium Weight Cargo (80\% ratio).
    \end{tablenotes}
\end{threeparttable}
\end{table}

The comparison with the stability computer shows that, in general, our implementation and simplification of vessel stability are effective. However, stability errors persist for all three vessels. The errors related to propeller immersion are due to the minimum draft specification for each vessel being set at the boundary of feasibility. Propeller immersion could be improved by increasing the draft limit or allowing the trim to change. Further, our model generated stowage plans under a zero-trim assumption, meaning it cannot trim the vessel to fully immerse the propeller in the water when cargo weight is low. The errors related to heel, GM, and KG are mainly due to low draft limits. Firstly, we do not account for heel directly, only indirectly from the TCG. Second, we do not include any free-surface effects which does impact the GM limits sometimes. Overall, the hybrid-slot model produces stowage plans that are generally stable and in line with industry standards, with most deviations attributable to parameter settings rather than fundamental modelling errors. \autoref{tab:detailed_stability_errors_rerun} presents the updated error table where the minimum draft is increased to show that we can eliminate most errors by using a more realistic draft.

\begin{figure}[h!]
    \centering
    \includegraphics[width=0.5\linewidth]{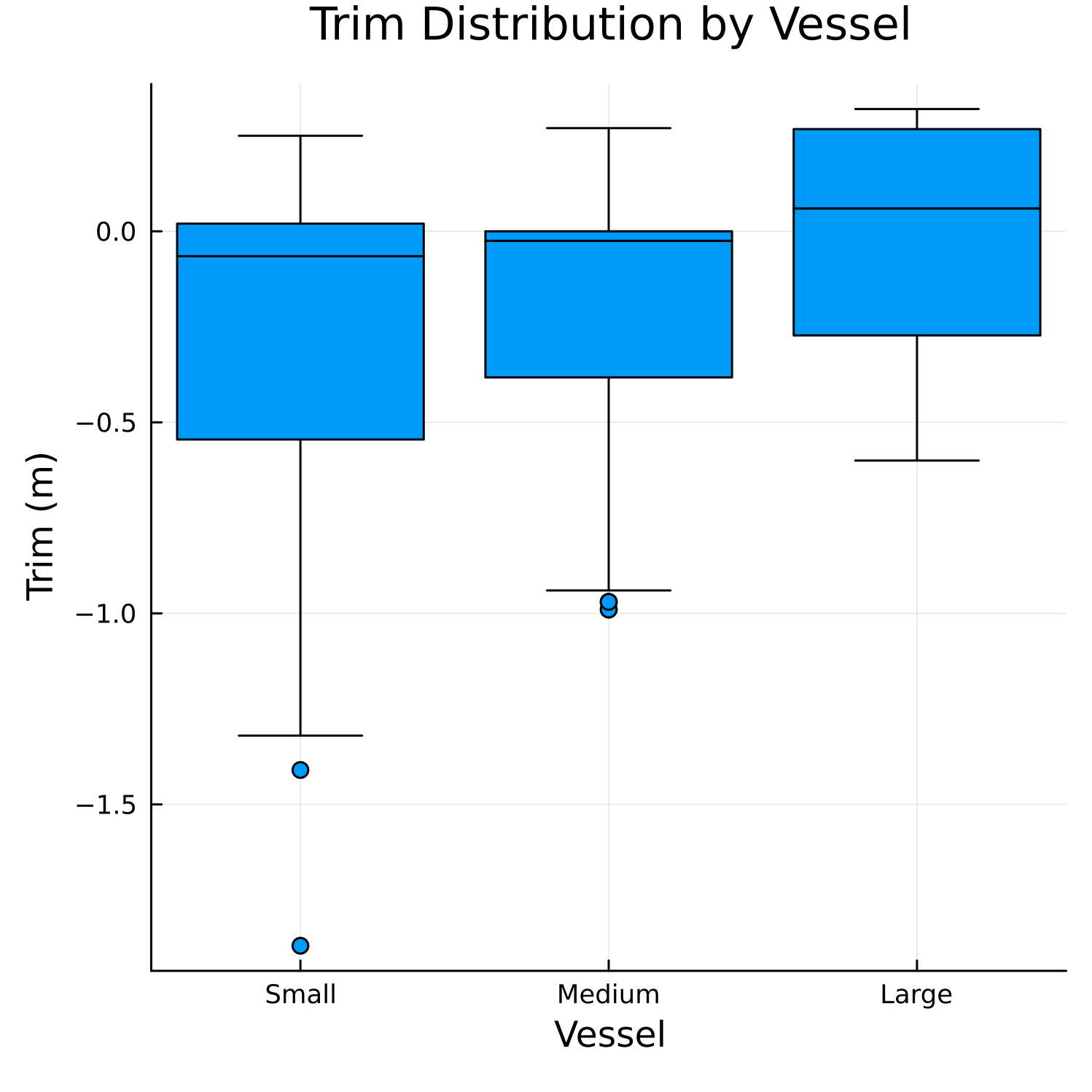}
    \caption{The boxplot show the trim distribution accross all solutions from the hybrid-slot model. The data is grouped by the vessel size. Positive trim is to the aft and negative to the front of the vessel.}
    \label{fig:trim_boxplot}
\end{figure}

Finally, we examine the validity of the zero-trim assumption in the hybrid-slot model. \autoref{fig:trim_boxplot} shows that the average across all instances and vessels is very close to zero. The model's trim generally biases towards the front of the vessel across all vessel sizes. However, the small and medium vessel exhibit high variance in trim between instances. This indicates that even though our calculations are quite precise, cumulative errors occur. To mitigate this, either the resolution of the hydrostatic data could be increased, the ballast tank centre of gravity could be approximated piecewise linearly, or the zero-trim assumption could be removed. However, the trim results in only 4 errors in the medium vessel. A key improvement to the stability formulations would be to introduce trim adjustments. This would enable greater freedom in the solutions generated while also aligning with industry practice. 

From an industry perspective, seaworthiness is important, but it must also be practical. In a real-world setting, the stability formulations presented here are likely too granular, as ballast water is rarely adjusted to the level of precision demonstrated here. This is due to operational constraints, including mechanical limitations on pump flow rates, inefficiencies in shifting negligible volumes of water, and regulatory requirements to treat ballast water prior to discharge. However, the precision of the formulation and implementation demonstrates that our assumptions regarding seaworthiness are valid. Furthermore, operational constraints can be incorporated, such as limiting the rate of ballast exchange or fixing specific ballast water levels based on the vessel's current state.  

\section{Conclusion} \label{sec:conclusion}
This study presented a comparative analysis of lane-based, grid-based, and hybrid-slot representations for RoRo stowage planning. We evaluated the three models with no stability, simple stability and full seaworthiness constraints across different vessel sizes and cargo configurations. The results show that the hybrid-slot model performs well across all the scenarios, but still does not solve everything to optimality within the time limit of an hour. The grid model performs second best but is limited by its computational complexity as the grid size increases. The lane-based model performs worst and finds optimal solutions in only a few cases for the small ship. This limitation is attributed to the combinatorial complexity of the precedence constraints. We further extended the hybrid-slot model to stow hazardous and refrigerated cargo using lazy constraints. Our results show that the implementation works well for realistic instances. Finally, we test the stability of our stowage plans from the hybrid-slot model against a loading computer. The results show that our stowage plans, in general, are within acceptable limits for industrial usage. Beyond the core comparative study, we extended the hybrid-slot model to accommodate hazardous and refrigerated cargo by implementing lazy constraints. Our results confirm that this approach is practical for realistic instances.
Furthermore, validation against a proprietary stability computer showed that the stowage plans generated by the hybrid-slot model generally remain within acceptable stability limits for industrial deployment, validating our underlying assumptions. Although the results are positive, computational times are not consistently viable for industrial deployment. For future work, we advise research on faster solution approaches such as heuristic methods. Additionally, attention must be paid to the robustness of these solutions. Mathematical optimisation tends to converge on the boundaries of feasibility; therefore, future work should quantify the operational flexibility of such plans, particularly regarding practical implementation under uncertain weather conditions and dynamic cargo bookings.

\section*{Conflicts of interest}
The authors declare no conflict of interest.

\newpage

\bibliographystyle{apalike-ejor}
\bibliography{references}

\section{Appendix}

\begin{table}[!htb]
\centering
\scriptsize 
\setlength{\tabcolsep}{2.5pt} 
\renewcommand{\arraystretch}{1.2}

\begin{threeparttable}
    \caption{Detailed performance breakdown of the hybrid-slot model across all stability configurations. 
    \textnormal{Gap = Optimality Gap (\%); Time in seconds; O/S = Optimal/Solved instances. Objective values omitted for brevity.}}
    \label{tab:detailed_breakdown_hybrid}

    \begin{tabular}{@{}ll cccc cccc cccc@{}}
    \toprule
    & & \multicolumn{4}{c}{\textbf{No Stability}} & \multicolumn{4}{c}{\textbf{Simple Stability}} & \multicolumn{4}{c}{\textbf{Full Stability}} \\
    \cmidrule(lr){3-6} \cmidrule(lr){7-10} \cmidrule(lr){11-14}
    \textbf{Ship} & \textbf{Inst.} & \textbf{Car\%} & \textbf{Gap\%} & \textbf{Time} & \textbf{O/S} & \textbf{Car\%} & \textbf{Gap\%} & \textbf{Time} & \textbf{O/S} & \textbf{Car\%} & \textbf{Gap\%} & \textbf{Time} & \textbf{O/S} \\
    \midrule

    \multirow{7}{*}{\textbf{Small}}  
         & \textbf{Avg.}  & \textbf{97.1}  & \textbf{0.02}  & \textbf{325.4} & \textbf{50/54} & \textbf{97.1}  & \textbf{0.16}  & \textbf{1001.2} & \textbf{41/54} & \textbf{97.1}  & \textbf{0.08}  & \textbf{944.6}  & \textbf{43/54} \\
         & ED             & 98.7           & 0.10           & 1324.5         & 6/9            & 98.6           & 0.62           & 1794.4          & 5/9            & 98.7           & 0.11           & 1670.4          & 6/9            \\
         & ML             & 98.2           & 0.03           & 602.3          & 8/9            & 98.2           & 0.10           & 2085.8          & 5/9            & 98.2           & 0.11           & 2005.5          & 5/9            \\
         & NC             & 96.8           & 0.00           & 8.2            & 9/9            & 96.8           & 0.21           & 1202.0          & 6/9            & 96.8           & 0.20           & 1202.8          & 6/9            \\
         & TO             & 94.9           & 0.00           & 0.2            & 9/9            & 94.9           & 0.00           & 1.0             & 9/9            & 94.9           & 0.00           & 1.1             & 9/9            \\
         & TS             & 96.0           & 0.00           & 0.4            & 9/9            & 96.0           & 0.03           & 459.1           & 8/9            & 96.0           & 0.00           & 319.3           & 9/9            \\
         & TH             & 98.0           & 0.00           & 16.9           & 9/9            & 98.0           & 0.03           & 465.2           & 8/9            & 98.0           & 0.03           & 468.3           & 8/9            \\
    \addlinespace

    \multirow{7}{*}{\textbf{Medium}} 
         & \textbf{Avg.}  & \textbf{95.7}  & \textbf{0.33}  & \textbf{832.6} & \textbf{45/54} & \textbf{95.1}  & \textbf{0.89}  & \textbf{1420.8} & \textbf{33/53} & \textbf{95.1}  & \textbf{0.79}  & \textbf{1381.0} & \textbf{34/53} \\
         & ED             & 98.3           & 1.49           & 2121.4         & 6/9            & 97.2           & 2.62           & 2451.8          & 3/9            & 97.0           & 2.75           & 2461.5          & 3/9            \\
         & ML             & 97.8           & 0.49           & 2420.8         & 3/9            & 96.7           & 2.50           & 2461.8          & 3/9            & 96.8           & 1.80           & 2493.3          & 3/9            \\
         & NC             & 96.0           & 0.00           & 182.0          & 9/9            & 96.0           & 0.06           & 886.5           & 7/9            & 96.0           & 0.05           & 884.1           & 7/9            \\
         & TO             & 90.9           & 0.00           & 0.7            & 9/9            & 89.7           & 0.04           & 801.6           & 6/8            & 89.7           & 0.04           & 802.6           & 6/8            \\
         & TS             & 94.0           & 0.00           & 2.9            & 9/9            & 94.0           & 0.06           & 825.6           & 7/9            & 94.0           & 0.06           & 821.9           & 7/9            \\
         & TH             & 97.0           & 0.00           & 268.0          & 9/9            & 96.9           & 0.07           & 1097.6          & 7/9            & 97.0           & 0.00           & 822.5           & 8/9            \\
    \addlinespace

    \multirow{7}{*}{\textbf{Large}}  
         & \textbf{Avg.}  & \textbf{96.6}  & \textbf{0.58}  & \textbf{854.9} & \textbf{39/48} & \textbf{94.6}  & \textbf{3.46}  & \textbf{1407.6} & \textbf{36/54} & \textbf{93.1}  & \textbf{6.18}  & \textbf{1514.7} & \textbf{33/53} \\
         & ED             & 97.3           & 2.84           & 2312.1         & 3/8            & 96.6           & 3.57           & 2663.2          & 3/9            & 89.2           & 17.17          & 2683.2          & 3/9            \\
         & ML             & 99.5           & 0.56           & 1006.3         & 3/4            & 88.5           & 16.64          & 2681.4          & 3/9            & 87.5           & 19.23          & 2824.5          & 3/9            \\
         & NC             & 96.7           & 0.00           & 300.5          & 9/9            & 96.7           & 0.03           & 933.9           & 8/9            & 96.7           & 0.06           & 1323.0          & 6/9            \\
         & TO             & 92.9           & 0.00           & 2.2            & 9/9            & 92.9           & 0.00           & 4.5             & 9/9            & 92.0           & 0.00           & 8.5             & 8/8            \\
         & TS             & 95.7           & 0.03           & 412.4          & 8/9            & 95.7           & 0.05           & 810.8           & 7/9            & 95.7           & 0.05           & 815.1           & 7/9            \\
         & TH             & 97.5           & 0.05           & 1096.0         & 7/9            & 97.4           & 0.50           & 1351.6          & 6/9            & 97.4           & 0.58           & 1433.8          & 6/9            \\ 

    \bottomrule
    \end{tabular}
\end{threeparttable}
\end{table}

\begin{table}[!htb]
\centering
\scriptsize 
\setlength{\tabcolsep}{2.5pt} 
\renewcommand{\arraystretch}{1.2}

\begin{threeparttable}
    \caption{Detailed performance breakdown of the grid model across all stability configurations. 
    \textnormal{Gap = Optimality Gap (\%); Time in seconds; O/S = Optimal/Solved instances. Objective values omitted for brevity.}}
    \label{tab:detailed_breakdown_grid}

    \begin{tabular}{@{}ll cccc cccc cccc@{}}
    \toprule
    & & \multicolumn{4}{c}{\textbf{No Stability}} & \multicolumn{4}{c}{\textbf{Simple Stability}} & \multicolumn{4}{c}{\textbf{Full Stability}} \\
    \cmidrule(lr){3-6} \cmidrule(lr){7-10} \cmidrule(lr){11-14}
    \textbf{Ship} & \textbf{Inst.} & \textbf{Car\%} & \textbf{Gap\%} & \textbf{Time} & \textbf{O/S} & \textbf{Car\%} & \textbf{Gap\%} & \textbf{Time} & \textbf{O/S} & \textbf{Car\%} & \textbf{Gap\%} & \textbf{Time} & \textbf{O/S} \\
    \midrule

    \multirow{7}{*}{\textbf{Small}}  
         & \textbf{Avg.}  & \textbf{85.8}  & \textbf{0.00}  & \textbf{9.5}   & \textbf{54/54} & \textbf{85.2}  & \textbf{0.02}  & \textbf{529.5}  & \textbf{50/54} & \textbf{84.5}  & \textbf{0.02}  & \textbf{569.2}  & \textbf{50/54} \\
                                 & ED             & 85.2           & 0.00           & 16.0           & 9/9            & 84.3           & 0.04           & 1317.3          & 7/9            & 83.3           & 0.04           & 1299.2          & 7/9            \\
                                 & ML             & 82.8           & 0.00           & 17.4           & 9/9            & 79.9           & 0.02           & 942.3           & 8/9            & 76.9           & 0.06           & 1669.6          & 7/9            \\
                                 & NC             & 89.1           & 0.00           & 9.9            & 9/9            & 89.1           & 0.03           & 671.5           & 8/9            & 89.1           & 0.00           & 217.5           & 9/9            \\
                                 & TO             & 84.1           & 0.00           & 1.4            & 9/9            & 84.1           & 0.00           & 5.0             & 9/9            & 84.1           & 0.00           & 8.2             & 9/9            \\
                                 & TS             & 86.3           & 0.00           & 2.2            & 9/9            & 86.3           & 0.00           & 8.1             & 9/9            & 86.3           & 0.00           & 12.9            & 9/9            \\
                                 & TH             & 87.5           & 0.00           & 10.3           & 9/9            & 87.5           & 0.00           & 232.7           & 9/9            & 87.5           & 0.00           & 207.8           & 9/9            \\
    \addlinespace

    \multirow{7}{*}{\textbf{Medium}} 
& \textbf{Avg.}  & \textbf{90.4}  & \textbf{0.00}  & \textbf{69.5}  & \textbf{54/54} & \textbf{89.2}  & \textbf{0.07}  & \textbf{1938.6} & \textbf{25/53} & \textbf{89.0}  & \textbf{0.07}  & \textbf{1950.2} & \textbf{25/52} \\
                                 & ED             & 89.5           & 0.00           & 106.5          & 9/9            & 88.4           & 0.10           & 3457.1          & 1/9            & 88.2           & 0.10           & 3482.2          & 1/9            \\
                                 & ML             & 87.5           & 0.00           & 90.8           & 9/9            & 82.3           & 0.09           & 3227.3          & 1/9            & 81.7           & 0.08           & 3189.1          & 1/8            \\
                                 & NC             & 93.2           & 0.00           & 113.3          & 9/9            & 93.2           & 0.14           & 2424.0          & 3/9            & 93.2           & 0.15           & 2431.6          & 3/9            \\
                                 & TO             & 89.7           & 0.00           & 12.8           & 9/9            & 88.4           & 0.00           & 38.9            & 8/8            & 88.4           & 0.00           & 65.1            & 8/8            \\
                                 & TS             & 91.2           & 0.00           & 24.9           & 9/9            & 91.2           & 0.00           & 55.6            & 9/9            & 91.2           & 0.00           & 90.6            & 9/9            \\
                                 & TH             & 91.6           & 0.00           & 68.5           & 9/9            & 91.6           & 0.10           & 2428.9          & 3/9            & 91.6           & 0.10           & 2442.6          & 3/9            \\
    \addlinespace

    \multirow{7}{*}{\textbf{Large}}  
& \textbf{Avg.}  & \textbf{91.0}  & \textbf{0.02}  & \textbf{555.0} & \textbf{50/54} & \textbf{91.8}  & \textbf{0.03}  & \textbf{1733.7} & \textbf{30/49} & \textbf{90.0}  & \textbf{0.03}  & \textbf{1882.6} & \textbf{29/53} \\
                                 & ED             & 89.7           & 0.00           & 407.6          & 9/9            & 90.6           & 0.08           & 3602.5          & 0/7            & 88.6           & 0.08           & 3602.9          & 0/9            \\
                                 & ML             & 87.9           & 0.00           & 295.7          & 9/9            & 91.4           & 0.02           & 2120.1          & 3/6            & 83.8           & 0.02           & 2760.0          & 3/9            \\
                                 & NC             & 94.3           & 0.10           & 2247.9         & 5/9            & 94.4           & 0.08           & 2490.5          & 3/9            & 94.4           & 0.06           & 2514.1          & 3/9            \\
                                 & TO             & 90.1           & 0.00           & 41.0           & 9/9            & 90.1           & 0.00           & 241.9           & 9/9            & 88.9           & 0.00           & 294.9           & 8/8            \\
                                 & TS             & 91.3           & 0.00           & 121.6          & 9/9            & 91.3           & 0.00           & 249.2           & 9/9            & 91.3           & 0.00           & 321.0           & 9/9            \\
                                 & TH             & 92.8           & 0.00           & 216.3          & 9/9            & 92.8           & 0.04           & 1698.0          & 6/9            & 92.8           & 0.04           & 1802.9          & 6/9            \\

    \bottomrule
    \end{tabular}
\end{threeparttable}
\end{table}

\begin{table}[!htb]
\centering
\scriptsize 
\setlength{\tabcolsep}{2.5pt} 
\renewcommand{\arraystretch}{1.2}

\begin{threeparttable}
    \caption{Detailed performance breakdown of the lane model across all stability configurations. 
    \textnormal{Gap = Optimality Gap (\%); Time in seconds; O/S = Optimal/Solved instances. Objective values omitted for brevity.}}
    \label{tab:detailed_breakdown_lane}

    \begin{tabular}{@{}ll cccc cccc cccc@{}}
    \toprule
    & & \multicolumn{4}{c}{\textbf{No Stability}} & \multicolumn{4}{c}{\textbf{Simple Stability}} & \multicolumn{4}{c}{\textbf{Full Stability}} \\
    \cmidrule(lr){3-6} \cmidrule(lr){7-10} \cmidrule(lr){11-14}
    \textbf{Ship} & \textbf{Inst.} & \textbf{Car\%} & \textbf{Gap\%} & \textbf{Time} & \textbf{O/S} & \textbf{Car\%} & \textbf{Gap\%} & \textbf{Time} & \textbf{O/S} & \textbf{Car\%} & \textbf{Gap\%} & \textbf{Time} & \textbf{O/S} \\
    \midrule

    \multirow{7}{*}{\textbf{Small}}  
& \textbf{Avg.}  & \textbf{69.3}  & \textbf{180.83} & \textbf{2836.2} & \textbf{12/51} & \textbf{68.7}  & \textbf{121.11} & \textbf{2906.3} & \textbf{11/43} & \textbf{91.6}  & \textbf{16.87} & \textbf{3404.9} & \textbf{7/12} \\
                                 & ED             & 33.1           & 360.11          & 3600.0          & 0/9            & 38.1           & 294.58          & 3600.0          & 0/3            & ---            & ---            & 3625.7          & 0/0           \\
                                 & ML             & 23.3           & 671.24          & 3600.0          & 0/6            & 25.7           & 355.33          & 3600.0          & 0/4            & ---            & ---            & 3622.8          & 0/0           \\
                                 & NC             & 91.2           & 10.60           & 2414.9          & 3/9            & 89.3           & 14.15           & 2445.2          & 3/9            & 100.0          & 0.00           & 2992.7          & 3/3           \\
                                 & TO             & 93.8           & 7.44            & 2409.6          & 3/9            & 92.5           & 9.54            & 2415.6          & 3/9            & 80.2           & 50.03          & 3253.1          & 2/3           \\
                                 & TS             & 90.9           & 11.08           & 2407.9          & 3/9            & 89.2           & 13.79           & 2432.4          & 3/9            & 88.4           & 15.00          & 3326.2          & 2/5           \\
                                 & TH             & 83.4           & 24.53           & 2564.3          & 3/9            & 77.5           & 39.29           & 2920.5          & 2/9            & 97.6           & 2.46           & 3609.0          & 0/1           \\
    \addlinespace

    \multirow{7}{*}{\textbf{Medium}} 
& \textbf{Avg.}  & \textbf{64.4}  & \textbf{107.17} & \textbf{3601.9} & \textbf{0/28}  & \textbf{61.2}  & \textbf{78.06}  & \textbf{3453.5} & \textbf{0/8}   & \textbf{0.0}   & \textbf{0.00}  & \textbf{3369.9} & \textbf{0/0}  \\
                                 & ED             & ---            & ---             & ---             & 0/0            & ---            & ---             & ---             & 0/0            & ---            & ---            & ---             & 0/0           \\
                                 & ML             & ---            & ---             & ---             & 0/0            & ---            & ---             & ---             & 0/0            & ---            & ---            & ---             & 0/0           \\
                                 & NC             & 52.0           & 196.67          & 3600.0          & 0/9            & 66.9           & 74.14           & 3600.0          & 0/3            & ---            & ---            & ---             & 0/0           \\
                                 & TO             & 70.2           & 124.47          & 3600.0          & 0/7            & 59.1           & 71.97           & 2892.9          & 0/3            & ---            & ---            & ---             & 0/0           \\
                                 & TS             & 68.7           & 53.00           & 3600.0          & 0/9            & 57.5           & 88.08           & 3600.0          & 0/2            & ---            & ---            & ---             & 0/0           \\
                                 & TH             & 66.6           & 54.56           & 3600.0          & 0/3            & ---            & ---             & ---             & 0/0            & ---            & 0.00           & 0.0             & 0/0           \\
    \addlinespace

    \multirow{7}{*}{\textbf{Large}}  
& \textbf{Avg.}  & \textbf{---}   & \textbf{---}    & \textbf{---}    & \textbf{0/0}   & \textbf{---}   & \textbf{---}    & \textbf{---}    & \textbf{0/0}   & \textbf{---}   & \textbf{---}   & \textbf{---}    & \textbf{0/0}  \\
                                 & ED             & ---            & ---             & ---             & 0/0            & ---            & ---             & ---             & 0/0            & ---            & ---            & ---             & 0/0           \\
                                 & ML             & ---            & ---             & ---             & 0/0            & ---            & ---             & ---             & 0/0            & ---            & ---            & ---             & 0/0           \\
                                 & NC             & ---            & ---             & ---             & 0/0            & ---            & ---             & ---             & 0/0            & ---            & ---            & ---             & 0/0           \\
                                 & TO             & ---            & ---             & ---             & 0/0            & ---            & ---             & ---             & 0/0            & ---            & ---            & ---             & 0/0           \\
                                 & TS             & ---            & ---             & ---             & 0/0            & ---            & ---             & ---             & 0/0            & ---            & ---            & ---             & 0/0           \\
                                 & TH             & ---            & ---             & ---             & 0/0            & ---            & ---             & ---             & 0/0            & ---            & ---            & ---             & 0/0           \\ \hline

    \bottomrule
    \end{tabular}
\end{threeparttable}
\end{table}

\begin{table}[!htb]
\centering
\small
\begin{threeparttable}
    \caption{Results of the Grid Model with No Stability (0.5$\times$0.5\,m resolution). 
    \textnormal{\footnotesize Objective function value is $\times 10^4$.}}
    \label{tab:grid_low_res}

    \begin{tabular}{@{}ll rrr r c@{}}
    \toprule
    \textbf{Ship} & \textbf{Instances} & \textbf{Cargo} & \textbf{Gap} & \textbf{Obj.} & \textbf{Time} & \textbf{Status} \\
     & & \textbf{(\%)} & \textbf{(\%)} & & \textbf{(s)} & \textbf{(Opt/Sol)} \\
    \midrule

    \multirow{7}{*}{\textbf{Small}} 
    & \textbf{Average} & \textbf{98.54} & \textbf{0.26} & \textbf{181.41} & \textbf{924.76} & \textbf{44/54} \\
    & ED               & 98.63          & 0.73          & 286.22          & 1485.95         & 6/9 \\
    & ML               & 99.05          & 0.53          & 318.00          & 1375.38         & 6/9 \\
    & NC               & 98.82          & 0.25          & 115.33          & 1312.94         & 6/9 \\
    & TO               & 96.32          & 0.00          & 108.00          & 122.97          & 9/9 \\
    & TS               & 98.66          & 0.00          & 102.00          & 285.86          & 9/9 \\
    & TH               & 99.77          & 0.06          & 158.89          & 965.49          & 8/9 \\

    \bottomrule
    \end{tabular}
    \begin{tablenotes}
        \footnotesize
        \item \textbf{Note:} Results are only present for the Small vessel. Medium and Large vessel instances exceeded the 40GB RAM limit during solving.
    \end{tablenotes}
\end{threeparttable}
\end{table}

\begin{table}[!htb]
\centering
\scriptsize 
\setlength{\tabcolsep}{3pt} 
\renewcommand{\arraystretch}{1.1}

\begin{threeparttable}
    \caption{Stability Errors by Cargo Configuration with improved draft limits. 
    \textnormal{View of instances violating stability limits.}}
    \label{tab:detailed_stability_errors_rerun}

    \begin{tabular}{@{}ll ll r c r@{}}
    \toprule
    \textbf{Vessel} & \textbf{Config} & \textbf{Distribution} & \textbf{Measure} & \textbf{Err.} & \textbf{Mean} & \textbf{Limit} \\
    \midrule

    \multirow{4}{*}{\textbf{Small}} 
    & \multirow{2}{*}{\textbf{Heavy 100}} 
      & No Cars & GM Actual & 1 & 1.16\,m & --- \\
    & & Trailer + Secu & GM Actual & 1 & 1.24\,m & --- \\
    \cmidrule{2-7}

    & \multirow{1}{*}{\textbf{Heavy 120}} 
      & \multirow{1}{*}{All Equal} & GM Actual & 1 & 1.21\,m & --- \\
    \midrule

    \multirow{1}{*}{\textbf{Medium}} 
    & \multirow{1}{*}{\textbf{Heavy 100}} 
      & \multirow{1}{*}{All Equal} & GM Actual & 1 & 1.01\,m & --- \\
    
    \bottomrule
    \end{tabular}
    
    \begin{tablenotes}
        \scriptsize
        \item \textbf{Notes:} Prop. Imm. = Propeller Immersion.
        \item Med. 80 = Medium Weight Cargo (80\% ratio).
    \end{tablenotes}
\end{threeparttable}
\end{table}

\end{document}